\documentclass[11pt]{article}
\usepackage{amsmath}  % Nodig voor het gebruik van \numberwithin
\usepackage{amsthm} % nodig voor het gebruik van \begin{proof}
\usepackage{amssymb} % nodig voor \Bbb en \mathbf (wat beide hetzelfde is, maar \Bbb wordt afgeraden)

\usepackage{hyperref} % nodig voor het commando \href{} gebruikt in de bibliografie

\usepackage{enumitem} % Dit pakketje laat toe om met itemize wat andere dingen te doen - zie andere artikels

\usepackage{color} % Voor het gebruik van kleuren

\newtheorem{thm}{Theorem}
\newtheorem{inspr}[thm]{}

\newenvironment{definitie}{\begin{itemize}\item[ ]\hspace{-26pt}\bf Definition \rm }{\end{itemize}}
\newenvironment{notatie}{\begin{itemize}\item[ ]\hspace{-26pt}\bf Notation \rm }{\end{itemize}}
\newenvironment{voorbeeld}{\begin{itemize}\item[ ]\hspace{-26pt}\bf Example \rm }{\end{itemize}}
\newenvironment{stelling}{\begin{itemize}\item[ ]\hspace{-26pt}\bf Theorem \rm }{\end{itemize}}
\newenvironment{propositie}{\begin{itemize}\item[ ]\hspace{-26pt}\bf Proposition \rm }{\end{itemize}}
\newenvironment{lemma}{\begin{itemize}\item[ ]\hspace{-26pt}\bf Lemma \rm }{\end{itemize}}
\newenvironment{opmerking}{\begin{itemize}\item[ ]\hspace{-26pt}\bf Remark \rm }{\end{itemize}}
\newenvironment{voorwaarde}{\begin{itemize}\item[ ]\hspace{-26pt}\bf Condition \rm }{\end{itemize}}
\newenvironment{oefening}{\begin{itemize}\item[ ]\hspace{-26pt}\bf Exercise \rm }{\end{itemize}}

\newenvironment{probleem}{\begin{itemize}\item[ ]\hspace{-26pt}\bf Problem \rm }{\end{itemize}}

\newcommand{\defin}{\begin{inspr}\begin{definitie}}  %\def already defined
\newcommand{\edefin}{\end{definitie}\end{inspr}}

\newcommand{\notat}{\begin{inspr}\begin{notatie}}  %\not already defined
\newcommand{\enotat}{\end{notatie}\end{inspr}}

\newcommand{\voorb}{\begin{inspr}\begin{voorbeeld}}  %\not already defined
\newcommand{\evoorb}{\end{voorbeeld}\end{inspr}}

\newcommand{\stel}{\begin{inspr}\begin{stelling}}
\newcommand{\estel}{\end{stelling}\end{inspr}}

\newcommand{\prop}{\begin{inspr}\begin{propositie}}
\newcommand{\eprop}{\end{propositie}\end{inspr}}

\newcommand{\lem}{\begin{inspr}\begin{lemma}}
\newcommand{\elem}{\end{lemma}\end{inspr}}

\newcommand{\opm}{\begin{inspr}\begin{opmerking}}
\newcommand{\eopm}{\end{opmerking}\end{inspr}}

\newcommand{\voorw}{\begin{inspr}\begin{voorwaarde}}
\newcommand{\evoorw}{\end{voorwaarde}\end{inspr}}

\newcommand{\oef}{\begin{inspr}\begin{oefening}}
\newcommand{\eoef}{\end{oefening}\end{inspr}}

\newcommand{\prob}{\begin{inspr}\begin{probleem}}
\newcommand{\eprob}{\end{probleem}\end{inspr}}

\newcommand{\bew}{\vspace{-0.3cm}\begin{itemize}\item[ ] \bf Proof\rm: }
\newcommand{\ebew}{\hfill $\qed$ \end{itemize}}

\newcommand{\snl}{\vskip 7pt} % Het is noodzakelijk dat er voor de instructie vspace een lege lijn staat
\newcommand{\nl}{\vskip 12pt} % Kunnen we dit hierin opnemen?

\newcommand{\ot}{\otimes}

\newcommand{\inv}{^{-1}}

\newcommand{\tl}{\triangleleft}
\newcommand{\tr}{\triangleright}

\newcommand{\tussenen}{\qquad\quad\text{and}\qquad\quad}

\numberwithin{thm}{section}   % Zorgt ervoor dat de nummering bestaat uit ?.?
\numberwithin{equation}{section} % Zorgt ervoor dat de nummering bestaat uit ?.?

\begin{document}

\centerline{\bf \Large From Hopf algebras to topological quantum groups }
\centerline{\bf \large A short history, various aspects and some problems \rm $^{(*)}$ }
\vspace{13pt}
\centerline{\it Alfons Van Daele \rm $^{(+)}$}
\bigskip\bigskip

{\bf Abstract} 
\nl
Quantum groups have been studied within several areas of mathematics and mathematical physics. This has led to different approaches, each of them with their own  techniques and conventions.
\snl
Starting with Hopf algebras, where there is a general consensus, moving in the direction of topological quantum groups, where there is no such consensus, it is easy to get lost. Not only many difficulties have to be overcome, but also several choices must be made. The way this is done is often confusing. Some choices even turn out to be rather annoying. 
\snl
As an introductory lecture at the conference on `{\it Topological quantum groups and Hopf algebras}' in 2016, we have explained these `{\it choices, difficulties and annoyances}' encountered on the road from Hopf algebras to topological quantum groups.
\snl
 In these notes, we  discuss more aspects of the development of locally compact quantum groups. We not only explain some of these difficulties in greater detail, but we also give  background information about the different steps, combined with some historical comments.  
\snl 
 We start with finite quantum groups and continue with discrete quantum groups, compact quantum groups and algebraic quantum groups. Multiplicative unitaries are an important side track before we finally arrive at locally compact quantum groups.
\snl
Along the way, we also formulate some interesting remaining problems in the theory. 
\nl
\nl
Date: {\it 14 January 2019}
\nl
\vskip 1cm
\hrule
\medskip
\begin{itemize}
\item[$^{(*)}$] Notes written for the Proceedings of the conference {\it Topological quantum groups and Hopf algebras} (Warsaw, November 2016)
\item[$^{(+)}$] Department of Mathematics, University of Leuven, Celestijnenlaan 200B,
B-3001 Heverlee, Belgium. E-mail: \texttt{alfons.vandaele@kuleuven.be}
\end{itemize}
\newpage

% Inleiding

\section{\hspace{-17pt}. Introduction} \label{s:intro} % \input{artikel1.tex}

Hopf algebras find their origin in algebraic topology (\cite{Hopf}, 1941). The general theory was further developed by others, resulting in the standard reference books by Sweedler (\cite{Sweedler}, 1969) and Abe (\cite{Abe}, 1977). For a more recent treatment we refer to the book by Radford (\cite{Radford}, 2012). 
\snl
Quantum groups, as they were called by Drinfel'd, had a breakthrough in 1986 with his work presented at the International Mathematical Congress in Berkeley (\cite{Drinfeld}, 1986). Independently, results of the same type were obtained  by Jimbo (\cite{Jimbo}, 1987). 
\snl
Parallel developments took place with the attempts to generalize Pontryagin duality for abelian locally compact groups (\cite{Pontryagin}, 1939)  to the non-abelian case. After many years of research, with various partial, but important solutions, Kac and Vainerman on the one hand (\cite{VaK - 1974}, 1974) and Enock and Schwartz (\cite{E-S1}, 1975 and \cite{E-S2}, 1992) on the other hand, developed the theory of Kac algebras. These are objects admitting a dual of the same kind and their duality includes the Pontryagin duality of locally compact abelian groups, as well as the intermediate results of this kind obtained earlier for non-abelian locally compact groups.  The first example of a finite-dimensional Kac algebra was obtained by Kac and Paljutkin (\cite{KP - 1966}, 1966). 
\snl
Unfortunately, as it turned out by the the work of Woronowicz on compact quantum groups (\cite{W2}, 1987 and \cite{W3}, 1998), in particular with his construction of the quantum $SU_q(2)$ group (\cite{W1}, 1987), the theory of Kac algebras was found to be too restrictive. The new examples are not Kac algebras because they have an antipode whose square is not the identity map. 
\snl
In the seventies, when the Kac algebras where developed, it seemed quite natural to assume an antipode $S$ satisfying $S^2=\iota$. We use $\iota$ for the identity map. All the known examples at that time did have this property. Recall that the antipode is for a Hopf algebra what the inverse is for a group. So this condition was indeed very reasonable. On the other hand, before this, examples of Hopf algebras with an antipode not satisfying $S^2=\iota$ where already known. See e.g.\ \cite{Taft}, 1971. But the underlying algebras for these examples are not operator algebras. The quantum $SU_q(2)$ was the first example, constructed in the framework of operator algebras, with an antipode not satisfying $S^2=\iota$. Remark in passing that there was at that time very little mutual knowledge among the algebraists and the analysts working with this type of abstract harmonic analysis.
\snl
These new discoveries initiated new research. The goal now was to find a more general notion to replace the Kac algebras. It had to include the new examples. And preferably, the concept had to be  self-dual  in the sense that again, the dual of the new object could be constructed within the same category. The duality should generalize the duality of Kac algebras (and so of course also the earlier forms of duality generalizing Pontryagin's duality for locally compact abelian groups).
\snl
It took many years before this goal was achieved. First there is the work of Masuda and Nakagami (\cite{M-N}, 1994), continued in collaboration with Woronowicz (\cite{M-N-W}, 2003). It must be mentioned that the latter was available many years before it actually got published. In the mean time, a theory of locally compact quantum groups was also developed by Kustermans and Vaes (\cite{Ku-V1}, 1999 and \cite{Ku-V2}, 2000). In \cite{M-N} the underlying operator algebras are von Neumann algebras, whereas in \cite{M-N-W}, they are C$^*$-algebras. The same is true for the locally compact quantum groups as treated by Kustermans and Vaes. They studied locally compact quantum groups in the von Neumann algebra framework in \cite{Ku-V3} (2003). 
\snl
Nowadays, it is generally accepted that locally compact quantum groups as in \cite{Ku-V2} provide the best and strongest results. 
\snl
Remark that, from a motivational point of view, it was more natural to search for a notion of a locally compact quantum group within a C$^*$-algebra framework. After all C$^*$-algebras are considered as quantized locally compact spaces. However, as we will see in Section \ref{s:lcqgrps}, the von Neumann algebraic theory is easier.
\snl
We would also like to mention here that, to our knowledge, it was Kirchberg, at a conference in Kopenhagen (\cite{Ki}, 1992), who was the first to suggest a theory with an antipode $S$, not satisfying $S^2=\iota$, but with a polar decomposition.
\snl
It is somewhat remarkable and worthwhile to mention that before a satisfactory notion of a locally compact quantum group was obtained, there was a focus on examples. Highly sophisticated ones were mainly obtained by Woronowicz and some of his coworkers. The study of multiplicative unitaries, as developed mainly by Baaj and Skandalis (\cite{Ba-Sk}, 1993), provided a framework for these examples. The Haar measures for some of them were only considered later in unpublished work (\cite{VD6}, 2001). But it showed finally that these examples complied with the later developed notion.
\snl
Above is only a short account of the different intermediate steps leading up to the final notion of locally compact quantum groups. We will say a  more about this further in the paper. But the reader is also advised to compare it with the information given in other comprehensive works on the theory. 
\nl
\bf Choices, difficulties and annoyances. \rm  
\nl
As we see already with this limited introduction, many people from different schools were involved. On the way, many choices had to be made. Therefore, it should not come as a surprise that this had led to a variety of conventions. One of the reasons is  that the final notion came only later, after partial results had been found and examples of different kinds were constructed. This fact is without any doubt also responsible for these differences. 
\snl
The theory in itself is already quite complicated.  It comes with many technical difficulties. The problem with the different conventions is not only that it makes things even more difficult, it also causes annoying situations. 
\nl
\bf Expository character of the paper. \rm 
\nl
When we started to write these notes, in the first place, we had the idea of discussing these problems. After all, the notes were intended to expand the material presented at the conference in Warsaw where this topic was chosen  for our  introductory talk. However, in the process of writing, we found that for doing this, we needed to include various aspects of the development of locally compact quantum groups and some historical comments. This turned out, not only to be necessary, but also of some interest in its own right. 
\snl
 We strongly believe that this will be helpful mostly for young researchers who have already some knowledge about the operator algebraic approach to quantum groups but want to get a deeper understanding of the different steps leading to this theory, starting from Hopf algebras and locally compact groups. Moreover, by  discussing the passage from Hopf algebras to locally compact quantum groups, in the end, when we briefly treat locally compact quantum groups, we can illustrate some of the features there by similar ones in the simpler cases.
\snl
This paper is mostly expository and contains no new results.  We omit many details but we suggest to complete some of them as an exercise for the reader (cf.\ items marked as \emph{Exercise}). 
\snl
Here and there we also propose  possible research problems (cf.\ items marked as \emph{Problem}). One of them, perhaps the most interesting and challenging one, is to develop a theory of locally compact quantum groups where the existence of the Haar measure is not part of the axioms, but a result, following from another set of natural axioms. We include some ideas to attack this problem. 
\snl
Let us also refer to the book of Timmermann, \emph{An invitation to quantum groups and duality. From Hopf algebras to Multiplicative Unitaries and beyond} (\cite{Timmermann}, 2008). In a sense, it treats the same material, but from a different angle. Here we emphasize on other aspects. The reader who wants to go deeper into this theory is advised to read his work as well and compare the treatment of the various steps in his work with the way it is presented in this more concise note.
\nl
\bf Content of the paper.\rm
\nl 
In \emph{Section} \ref{s:fqgrps} we start our discussion with \emph{finite quantum groups} and we see that even here, some choices have to be agreed upon. In this paper, we will call a finite-dimensional Hopf $^*$-algebra a finite quantum group only if the underlying algebra is an operator algebra (see Definition \ref{defin:2.3} in Section \ref{s:fqgrps}).  We motivate this choice in detail, refer to the forthcoming more general case and use an example to make our point. Finite groups fit into this framework, but again, also here a choice has to be made.
\snl
In \emph{Section} \ref{s:dac} we proceed one step further in two different directions. We discuss the possible notions of discrete quantum groups and of compact quantum groups. For historical reasons, we start with looking for the good notion of a compact quantum group and arrive at the definition given by Woronowicz in \cite{W3}. We formulate some of the main results and use them to discuss other approaches to compact quantum groups. 
\snl
Discrete quantum groups were in the first place studied as duals of compact quantum groups, see \cite{P-W} (1990). Later the notion was defined without reference to a compact quantum group, see e.g. \cite{E-R} (1994) and \cite {VD3n} (1996). Thereby it was shown that the dual of a discrete quantum group is a compact quantum group as expected. It should be mentioned that both the compact and discrete quantum groups can be studied in a purely algebraic context, without reference to the operator algebras that they are part of. This leads to the more general notion of an algebraic quantum groups.
\snl
This is what we treat in \emph{Section} \ref{s:alg}. There is some confusion about the terminology here, but we will define an algebraic quantum group as a multiplier Hopf $^*$-algebra with positive integrals. Then the dual can be defined within the same category and the duality includes the duality between compact and discrete groups. Indeed, any compact quantum group is essentially a Hopf $^*$-algebra with positive integrals. On the other hand, it is shown that a discrete quantum group, as we define it in Section \ref{s:dac}, is a multiplier Hopf $^*$-algebra with positive integrals, admitting a cointegral.
\snl
Algebraic quantum groups are nice objects from different perspectives. First, it is a self-dual theory, including compact as well as discrete quantum groups. Secondly, it is a purely algebraic theory still with many (though not all) features in common with the general theory of locally compact quantum groups. Furthermore, it is far more easy to work with. We claim that first learning about algebraic quantum groups in the sense of this paper will facilitate to a great extend the later study of locally compact quantum groups. In fact, it is already useful to get a deeper understanding of discrete and compact quantum groups with their duality.
\snl
In this section, we also include a discussion of the Larson-Sweedler theorem. Roughly speaking it says that the antipode exists and can be obtained from the existence of integrals. This is a very important result in view of the theory of locally compact quantum groups as it is known today. We include here some related results that are important for understanding the development of locally compact quantum groups as we explain in Section \ref{s:lcqgrps}.
\snl
In the section on algebraic quantum groups, we also encounter the multiplicative unitaries. There are various approaches. This is treated in \emph{Section} \ref{s:mu}. It is possible to approach multiplicative unitaries from a purely algebraic point of view in the case of algebraic quantum groups. We start this section by doing so. However, multiplicative unitaries are mostly studied as unitary operators on Hilbert spaces. They appeared already a long time before the term itself was introduced by Baaj and Skandalis in \cite{Ba-Sk}. For this reason we spend some time to look at the history of this concept. Finally, we discuss the notion of a manageable multiplicative unitary as introduced by Woronowicz in \cite{W4}, see also \cite{S-W} and related concepts.
\snl
In Section \ref{s:lcqgrps} we arrive at the locally compact quantum groups, our goal for this paper. There are two operator algebraic approaches. On the one hand, there is a C$^*$-algebraic approach, while on the other hand, we have a von Neumann algebraic one. We begin the section with the two definitions and the relation between them. Then we focus on the von Neumann algebraic setting and we indicate how ideas from the purely algebraic theory of algebraic quantum groups, are used to develop the theory as it is done in \cite{VD9n}. For these ideas, we refer to Section \ref{s:alg} where we have included the relevant results for this purpose already.
\snl
In this section, we also compare the notion of a von Neumann algebraic locally compact quantum group with the earlier notion of a Kac algebra. And we formulate a problem about this also.
\snl
We finish the paper with some concluding remarks in Section \ref{s:conclusions}.  We also formulate a few suggestions about how to improve the theory and its presence in the literature. 
\nl
\bf Notations and conventions \rm
\nl
In general, our algebras are not required to have an identity, but we need that the product is non-degenerate.  
We use $1$  for the identity in an algebra when it exist and for the identity in its multiplier algebra if not. As mentioned already, we use $\iota$ to denote the identity map. The identity element in a group is denoted with $e$. Most of our algebras will be algebras over the field of complex numbers and have involutions.
\snl
In general we write $A\ot B$ for the tensor product of two algebras $A$ and $B$. In the case of operator algebras, we will use $A\odot B$ for the \emph{algebraic} tensor product while we keep using $A\ot B$ for the adapted \emph{completed} tensor product, depending on the context. 
\snl
We will use the leg-numbering notation in various circumstances. See e.g.\ Section \ref{s:mu}.
\snl
A coproduct on a finite-dimensional algebra $A$ is a homomorphism from $A$ to the tensor product $A\ot A$. When the algebra is no longer finite-dimensional we need an adapted notion. Also for an operator algebra this is the case. So, whenever we talk about a coproduct, it depends on the setting what is really meant.
\snl
For a coproduct we will use the Sweedler notation whenever this is convenient. This presents no problem in the case of a Hopf algebra. Using the Sweedler notation for multiplier Hopf algebras has been documented in the literature. See e.g.\ \cite{D-VD}, \cite{VD8n} and the forthcoming note \cite{VD12n}. In the case of an operator algebra however, it makes little sense to use it, except possibly for motivational reasons.
\snl
If $\mathcal H$ is a Hilbert space, we will use $\mathcal B(\mathcal H)$ for the von Neumann algebra of all bounded linear operators on $\mathcal H$. We denote the space of normal linear functionals on $\mathcal B(\mathcal H)$, the predual of this von Neumann algebra, by $\mathcal B(\mathcal H)_*$.
\nl
\bf Basic references \rm 
\nl 
For the theory of Hopf algebras we refer to the basic works of Abe \cite{Abe} and Sweedler \cite{Sweedler}, as well as to the more recent work by Radford \cite{Radford}. For compact quantum groups we have the work of Woronowicz \cite{W1, W2, W3} and for discrete quantum groups \cite{P-W} and \cite{VD3n}. Also the notes on compact quantum groups in \cite{M-VD} should be helpful. The basic reference for multiplier Hopf algebras is \cite{VD1} and for algebraic quantum groups (multiplier Hopf algebras with integrals), it is \cite{VD3}. For locally compact quantum groups, there is the work of Masuda, Nakagami and Woronowicz \cite{M-N,M-N-W} and the work of Kustermans and Vaes \cite{Ku-V1,Ku-V2,Ku-V3}.
\snl
For the general theory of C$^*$-algebras and von Neumann algebras, there are several possible references. One can look at the books of Kadison and Ringrose \cite{K-R} and at the books of Takesaki \cite{Tak1,Tak2}.
\nl
\bf Acknowledgments \rm
\nl
I would like to thank the organizers of the conference, {\it Topological quantum groups and Hopf algebras} (Warsaw, November 2016), for the opportunity to give the introductory lecture of the conference.

 \section{\hspace{-17pt}. Finite quantum groups } \label{s:fqgrps} %\input artikel2.tex 

What is a finite quantum group? Some people say that any finite-dimensional Hopf algebra is a finite quantum group. Other people will insist to have a finite-dimensional Hopf $^*$-algebra. In fact, then it is quite common to only admit operator algebras for the underlying $^*$-algebra. in this section, we will make a clear choice and motivate it.

\subsection{Finite-dimensional Hopf algebras}
First recall the following definition. Then consider an example to illustrate the upcoming choice.

\defin\label{defin:2.1}
A \emph{Hopf algebra} over a field $k$ is a pair $(A,\Delta)$ where $A$ is a unital $k$-algebra and $\Delta$ a coproduct on $A$ such that there is a counit $\varepsilon$ and an antipode $S$. When $A$ is a $^*$-algebra (over the field $\mathbb C$) and $\Delta$ is a $^*$-homomorphism, we call it a \emph{Hopf  $^*$-algebra}.
\edefin

For the definition of a Hopf algebra we refer to the standard works \cite{Abe, Sweedler} and also \cite{Radford}. Usually, the counit and the antipode are included explicitly in the definition, but because they are unique if they exist, this is not really necessary. We prefer to use the above formulation as this is closer to the upcoming definitions in the framework of operator algebras.
\snl
In the field of Hopf algebras, the study of Hopf $^*$-algebras is not very common, but it is of course very relevant when looking at the operator algebraic approach to quantum groups. 
\snl
Next consider the following example originally due to Taft \cite{Taft}.

\voorb\label{voorb:2.2}
Let $A$ be the unital algebra over $\mathbb C$ generated by elements $a,b$ satisfying  $a^4=1$, $ab=iba$ and $b^2=0$. A coproduct $\Delta$ on $A$ can be defined by
\begin{equation*}
\Delta(a)=a\ot a
\tussenen
\Delta(b)=a\ot b + b\ot a\inv
\end{equation*}
The pair $(A,\Delta)$ is a finite-dimensional Hopf algebra. 
It is a Hopf $^*$-algebra if we let $a$ and $b$ be self-adjoint.
\evoorb

The algebra is $8$-dimensional. The counit satisfies $\varepsilon(a)=1$ and $\varepsilon(b)=0$. For the antipode we have $S(a)=a\inv$ and $S(b)=ib$. Also $S^2(b)=-b$ so that for this example, we do not have that $S^2=\iota$. Instead we have $S^4=\iota$. Also observe that $S(b)^*=-ib$ so that $S(S(b)^*)^*=b$ as it should for a Hopf $^*$-algebra.
\snl
Consider the subalgebra $B$ of $A$ generated by $a^2$ and $ab$. Put $p=a^2$ and $q=ab$. Then $p^2=1$, still $q^2=0$ and  now $pq=-qp$. For the coproduct we find
\begin{equation*}
\Delta(p)=p\ot p
\tussenen
\Delta(q)=p\ot q + q\ot 1
\end{equation*}
and we see that $\Delta(B)\subseteq B\ot B$. Still $\varepsilon(p)=1$ and $\varepsilon(q)=0$ while $S(p)=p$ and $S(q)=iba\inv=-pq$. Remark that $q^*=ba=-iab=-iq$. The pair $(B,\Delta)$ is a $4$-dimensional Hopf $^*$-algebra. Again we have $S^4=\iota$ but not $S^2=\iota$ because $S^2(q)=S(pq)=pqp=-q$.
\snl
Remark that we can endow this algebra with a different $^*$-structure by letting $q^*=q$ as wel as $p^*=p$. This comes down to multiplying $q$ with a suitable number of modulus $1$. It is not essentially different. 
\snl
The  $^*$-algebra $A$ is not an operator algebra because we have a non-zero self-adjoint element $b$ with $b^2=0$. Also the subalgebra $B$ is not an operator algebra, essentially for the same reason. As a matter of fact, a four dimensional operator algebra has to be abelian or equal to the algebra of two by two matrices. And in the latter case, we can not have a  counit. Later we will see that $S^2=\iota$ is necessarily true if we have a finite-dimensional Hopf $^*$-algebra with an operator algebra, see Proposition \ref{prop:2.11}.

\oef
To get familiar with finite quantum groups it is instructive to take some time to verify the above statements. Look also at the related examples in \cite{VD4n}.
\eoef

We see from these examples that, if we insist to have operator algebras, then not all finite-dimensional Hopf $^*$-algebras will be finite quantum groups.
\snl
Here, we agree on the following more restricted definition.

\defin\label{defin:2.3}
Let $(A,\Delta)$ be a finite-dimensional Hopf $^*$-algebra. We call it a \emph{finite quantum group} if the underlying $^*$-algebra is an operator algebra.
\edefin

This means that $A$ has to be a direct sum of matrix algebras.

\opm
This is a \emph{first choice} we have to make. It is a natural one. If we agree that a finite group is also a locally compact group, we want that a finite quantum group is also a locally compact quantum group. Then we need to assume that the underlying algebra is an operator algebra.
\eopm

This is \emph{not the end of the story} about finite quantum groups as we will explain next.

\subsection{Finite groups and finite quantum groups}
Clearly, we want any finite group to be a finite quantum group. There are however, two ways to achieve this as we conclude from the two following well-known results.

\prop\label{prop:2.4}
Let $G$ be a finite group. Consider the $^*$-algebra $C(G)$ of all complex functions on $G$ with pointwise operations. The product in $G$ gives rise to a coproduct $\Delta$ on $C(G)$ defined by $\Delta(f)(p,q)=f(pq)$ where $p,q\in G$. The pair  $(C(G),\Delta)$ is a finite quantum group.
\eprop

\prop\label{prop:2.5}
Let $G$ be a finite group. Consider the group algebra $\mathbb C G$ and denote by $p\mapsto \lambda_p$ the canonical embedding of $G$ in $\mathbb C G$. It is a $^*$-algebra with $\lambda_p^*=\lambda_{p\inv}$. There is a coproduct $\Delta$ on $\mathbb C G$ defined by $\Delta(\lambda_p)=\lambda_p\ot\lambda_p$ for all $p\in G$. The pair $(\mathbb C G,\Delta)$ is again a finite quantum group.
\eprop 

Which one do we choose?  
\opm\label{opm:2.6}
i) In general, within the theory of locally compact quantum groups, it is most natural to choose the pair $(C( G),\Delta)$ as the finite quantum group associated to a finite group.\\
ii) This  is in agreement with the \emph{quantization procedure}: We start with a space together with certain properties. We take an (abelian) algebra of well-chosen functions on this space with the induced properties. And then we deform the algebra to become a non-abelian one, while we consider the properties on the deformed level. \\
iii) This is not always the choice made by the algebraists. And even some operator algebraists in the past have considered the other possibility or have not been consequent in their choices.
\eopm

\subsection{Integrals and cointegrals}
We can illustrate the above ambiguity using the notions of integrals and cointegrals.
\snl
The Haar measure on a finite group is the discrete measure and gives rise to the Haar integral $f\mapsto \sum_p f(p)$. It leads us to the notion of an \emph{integral on} a Hopf algebra.

\defin\label{defin:2.6}
Let $(A,\Delta)$ be a Hopf algebra. A non-zero linear functional $\varphi$ on $A$ is called a \emph{left integral} if $(\iota\ot\varphi)\Delta(a)=\varphi(a)1$ for all $a$. A non-zero linear functional $\psi$ on $A$ is called a \emph{right integral} if 
$(\psi\ot\iota)\Delta(a)=\psi(a)1$ for all $a$. Recall that we use $\iota$ for the identity map.
\edefin
There is also the notion of a \emph{cointegral in} a Hopf algebra.
\defin \label{defin:2.7}
Let $(A,\Delta)$ be a Hopf algebra with counit $\varepsilon$. A non-zero element $h$ in $A$ is called a \emph{left cointegral} if $ah=\varepsilon(a)h$ for all $a$ in $A$. A non-zero element $k$ in $A$ is called a \emph{right cointegral} if $ka=\varepsilon(a)k$ for all $a$.
\edefin

The algebraists use a different terminology than the one above. They will use an \emph{integral on} in the first case and an \emph{integral in} in the second case. This is clearly related with the ambiguity in the choice of the finite quantum group associated to a group as we discussed earlier. We claim that speaking about integrals and cointegrals as in Definitions \ref{defin:2.6} and \ref{defin:2.7} is better, given the choice of associating the pair $(C(G),\Delta)$ to the group $G$ as we have done.
\snl
In the case of a finite quantum group (as we defined it), we get the following result.

\prop\label{prop:2.11}
Let $(A,\Delta)$ be a finite quantum group. The antipode satisfies $S^2=\iota$. The left integral $\varphi$ is also right invariant and it is a trace.  We have  $\varphi(1)\neq 0$ and if we normalize it so that $\varphi(1)=1$, it is a positive linear functional.
\eprop

For a simple proof of this result, we refer to Section 2 of \cite{VD4n}.
\snl
In fact, if we have a finite-dimensional Hopf $^*$-algebra with a positive integral, then it is a finite quantum group. The reason is that integrals are automatically faithful and if there is a faithful positive linear functional, the $^*$-algebra is an operator algebra.

\voorb
Consider again the example with $p$ and $q$ and $p^2=1$, $q^2=0$ and $pq=-qp$. Define $\varphi$ on the algebra by $\varphi(1)=\varphi(p)=\varphi(q)=0$, while $\varphi(pq)=1$. This is a left integral. If we put $\psi(1)=\psi(p)=\psi(pq)=0$ while $\psi(q)=1$, we get the right integral. We see that $\varphi\neq\psi$ and that $\varphi(1)=\psi(1)=0$. This also shows that we do not have a finite quantum group  as defined in Definition \ref{defin:2.3}.
\evoorb

The integrals can not be positive as $\varphi(1)=\psi(1)=0$. The  integrals are not self-adjoint. However, this is not a big deal. Indeed, given e.g.\ a left integral $\varphi$, we will have that $\overline\varphi$, defined by $\overline\varphi(a)=\overline{\varphi(a^*)}$ for all $a$, will again be a left integral and so a scalar multiple of the original $\varphi$. Hence, by multiplying $\varphi$ with a well-chosen complex number, in this case with $i$, we get a self-adjoint integral.  
\snl
It also shows  that for this example we can not define another $^*$-structure, compatible with the coproduct, and so that the algebra is an operator algebra. 
\opm
The property $S^2=\iota$ shows that the finite-dimensional case is \emph{not a good model} for the general theory. It is too restrictive. In fact, in finite dimensions, it is better to look at the \emph{non-involutive} Hopf algebra's. The general locally compact quantum groups have more common features with those.
\eopm

An example as in Example \ref{voorb:2.2} can be used to illustrate this statement. Indeed, in infinite dimensions, it is possible to find locally compact quantum groups based on this type of quantization of the $ax+b$-group, see e.g. Remark \ref{rem:5.21}  further in Section \ref{s:mu}.  

%%%%%% Discrete and compact quantum groups %%%%%

 \section{\hspace{-17pt}. Discrete and compact quantum groups}  \label{s:dac} %\input artikel3.tex

The next questions are: What is a \emph{discrete quantum group}? What is a \emph{compact quantum group}? 
\snl
If we follow the spirit of quantization, a discrete quantum group should be a \emph{discrete quantum space} with a group-like structure while a compact quantum group is expected to be a \emph{compact quantum space}, again with an adapted group-like structure. 
\snl
In both cases, the quantum space should be an operator algebra of a certain kind and the group-like structure a coproduct on this algebra with some extra conditions (like the existence of a counit and an antipode). The notions should be such that, if a discrete quantum group is also a compact one, it should be a finite quantum group (in the sense of our Definition \ref{defin:2.3}).
\snl
For historical reasons we begin the discussion with compact quantum groups.  

\subsection{Compact quantum groups. Preliminary considerations}

We will see that already \emph{several difficulties} arise and again \emph{choices} have to be made. This is what we discuss first.
\snl
One aspect is easy. There is a wide consensus about the notion of a compact quantum space. It is a unital C$^*$-algebra. For the correct notion of a coproduct on a C$^*$-algebra, we get the inspiration from the group case.

\defin\label{defin:3.1}
Let $G$ be a compact group. Denote by $C(G)$ the C$^*$-algebra of continuous complex functions on $G$. We identify $C(G\times G)$ with the C$^*$-tensor product $A\ot A$. The product on $G$ gives rise to a coproduct $\Delta:A\to A\ot A$ by the formula $\Delta(f)(p,q)=f(pq)$ where $p,q\in G$.
\edefin

This naturally leads to the following notion.

\defin
A \emph{coproduct} on a unital C$^*$-algebra $A$ is a unital $^*$-homomorphism $\Delta:A\to A\ot A$ satisfying \emph{coassociativity} $(\Delta\ot\iota)\Delta=(\iota\ot\Delta)\Delta$.
\edefin

It is a common practice to use the minimal C$^*$-tensor product. This is motivated by the examples. As before $\iota$ is the identity map. The maps $\Delta\ot\iota$ and $\iota\ot\Delta$, in the first place defined on the algebraic tensor product $A\odot A$, stand here for the unique continuous extensions to maps from the minimal C$^*$-tensor product $A\ot A$ to the triple minimal tensor product $A\ot A\ot A$.

The next question to answer is then the following. What kind of conditions do we impose on the pair $(A,\Delta)$ to call it a compact quantum group? In particular, what does it mean here that the coproduct is \emph{group-like}?

It would be most natural to require the existence of a suitable counit and antipode. But this is less obvious than one might think. Indeed, here the \emph{difficulties already begin}.

\opm
i) In the motivating case of a compact group $G$ as in Definition \ref{defin:3.1} the unit $e$ in $G$ gives rise to a counit $\varepsilon:C(G)\to \mathbb C$ given by $\varepsilon(f)=f(e)$. This is a $^*$-homomorphism. However, from the examples we know, we can not expect the existence of a counit that  is an everywhere defined bounded $^*$-homomorphism from $A$ to $\mathbb C$. For the dual of $C(G)$, where the underlying algebra is the reduced group C$^*$-algebra $C_r^*(G)$, the obvious candidate for the counit is (in general) not everywhere defined.\\
ii) Because of this problem, we have to be careful with formulas like $(\iota\ot\varepsilon)\Delta(a)=a$ and $(\varepsilon\ot\iota)\Delta(a)=a$ for $a$ in $A$. It is not immediately clear what they mean because we have no standard procedure to extend the maps $\iota\ot\varepsilon$ and $\iota\ot\varepsilon$ to the completed tensor product $A\ot A$.\\
iii) Again in the group case, the inverse $p\mapsto p\inv$ gives rise to an antipode $S:C(G)\to C(G)$ given by $S(f)(p)=f(p\inv)$ for all $p\in G$. In this case, it is a $^*$-isomorphism.  This is even more problematic, for various reasons: 
\begin{itemize}
\item[-] From Hopf algebra theory, we know that we can not expect the antipode to be a homomorphism. It is an anti-homomorphism. This is not seen in the group case because the algebra $C(G)$ is abelian. There are also examples where the antipode is again not everywhere defined and unbounded. Now the dual of $C(G)$, the group algebra, is not an example but one can look at the quantum $SU_q(2)$ of Woronowicz. Remark that this problem was not existing for the earlier theory of Kac algebra. 
\item[-] The antipode is not expected to be a $^*$-map either but rather a bijective map that will satisfy $S(a)^*=S\inv(a^*)$ for elements $a$ in the domain of the antipode. This behavior with respect to the involution is rather natural as we know from Hopf $^*$-algebra theory.
\end{itemize}
iv) Because of these problems, there is no way to get a straightforward interpretation of the characterizing formulas
\begin{equation*}
m(S\ot \iota)\Delta(a)=\varepsilon(a)1
\tussenen
m(\iota\ot S)\Delta(a)=\varepsilon(a)1
\end{equation*}
for all $a$ in $A$. Here $m$ stands for the multiplication map.
\begin{itemize}
\item[-] The maps $S\ot\iota$ and $\iota\ot S$ would only be defined on a subspace of the algebraic tensor product $A\odot A$. They need not  be continuous and so there is no standard way to extend them to  maps on the completed tensor product.  
\item[-] The fact that the antipode is not a homomorphism, but an anti-homomorphism, causes extra problems.
\item[-] Finally, whereas in the case of an abelian C$^*$-algebra $A$, the multiplication map $m$ is well-defined as a bounded map from $A\ot A$ to $A$, this is no longer the case for non-abelian C$^*$-algebras. One can see this by looking at the multiplication map on the algebra $M_n$ of $n$ x $n$ matrices over $\mathbb C$. With increasing $n$, the norm of the multiplication map tends to infinity. 
\end{itemize}
One can say that the antipode is \emph{a fundamental}, if not \emph{the} problem for developing a theory of locally compact quantum groups. Solutions exist but are not obvious and sometimes rather annoying also.
\eopm

Indeed, the problems with the antipode caused \emph{great difficulties} in the development of a suitable notion of a locally compact quantum group. They are greatly responsible for the fact that it took so long before a good notion was found.

For Kac algebras, as developed in the 70s, independently by Kac and Vainerman on the one hand and Enock and Schwartz on the other hand, the following choices were made. We formulate them as part of the definition.

\defin (Incomplete)
A Kac algebra is a von Neumann algebra $M$ with a coproduct $\Delta:M\to M\ot M$.  The antipode, now called coinverse and denoted  by $\kappa$, is a $^*$-map with the property that $\kappa^2=\iota$. It is an anti-automorphism that flips the coproduct $\Delta$.
\edefin

For a precise definition of a Kac algebra we refer to Definition \ref{defin:6.2} in Section \ref{s:lcqgrps}.

The above is not the complete definition and so one should rather read it as a proposition. But it contains enough information for the discussion here. The tensor product $M\ot M$ considered is the von Neumann tensor product. We will comment later on the use of von Neumann algebras instead of C$^*$-algebras (see Remark \ref{opm:6.3} in Section \ref{s:lcqgrps}). Here we focus on the properties of the antipode. 

The following is an important remark.

\opm
The fact that the antipode is a $^*$-anti-automorphism that flips the coproduct is by no means sufficient for characterizing the antipode. If e.g.\ the underlying algebra is abelian and if moreover the coproduct is coabelian, the identity map will have the same properties.
\eopm

Within the theory of Kac algebras, this problem is overcome by a condition involving the integrals. But this is far from a natural requirement for an object that needs to replace the inverse in a group. Again see Definition \ref{defin:6.2}. This solution to the problem is \emph{not satisfactory}.

\subsection{Compact quantum groups}
In the middle of the 80s, Woronowicz first developed the theory of \emph{compact matrix pseudo groups} (\cite{W2}, 1987) and later \emph{compact quantum groups} (\cite{W3}, published in 1998 but developed earlier). Here is the definition.

\defin\label{defin:3.6}
Let $A$ be a \emph{unital} C$^*$-algebra and $\Delta$ a coproduct on $A$. The pair $(A,\Delta)$ is called a \emph{compact quantum group} if the sets $\Delta(A)(1\ot A)$ and $\Delta(A)(A\ot 1)$ are dense in $A\ot A$. 
\edefin

We use $\Delta(A)(1\ot A)$ for the linear span of elements of the form $\Delta(a)(1\ot b)$ where $a,b\in A$. Similarly for $\Delta(A)(A\ot 1)$. 

\opm
We observe the following:
\begin{itemize}[noitemsep]
\item[-] There is \emph{no counit} and \emph{no antipode} in the definition.
\item[-] The density conditions reflect dual forms of the cancellation law in a group.
\item[-] And basically, the underlying idea is that a compact semigroup with cancellation is a group. See \cite{Hofm}, also Proposition 3.2 in \cite{M-VD}. 
\end{itemize}
\eopm

Before continuing with the discussion, we formulate two of the main results in the theory of compact quantum groups as defined above.

\stel
Let $(A,\Delta)$ be a compact quantum group. There is a unique (non-zero) positive linear functional $\varphi$, satisfying left and right invariance, normalized so that $\varphi(1)=1$.
\estel

For the proof we refer to the original work by Woronowicz \cite{W3}. See also \cite{VD2n} and \cite{M-VD} where the proof is given for a slightly more general case. Remark that with this definition, the integral $\varphi$ on the C$^*$-algebra is not necessarily faithful.

There is nice theory of (co)-representations which is very similar to the representation theory of compact groups. The existence of the Haar state from the previous theorem is essential to obtain these results. We can consider polynomial functions defined as matrix elements of finite-dimensional representation. They form a $^*$-subalgebra that we denote as $\mathcal A$. It has the following properties.

\stel\label{thm:3.9}
The algebra $\mathcal A$ is a dense $^*$-subalgebra of $A$. It is invariant under the coproduct and the pair $(\mathcal A,\Delta)$ is a Hopf $^*$-algebra. It is the unique Hopf $^*$-algebra that is dense in $(A,\Delta)$.
\estel

The results about corepresentations lead us naturally to the dual discrete quantum group. We come back to this in a next subitem. First we focus on another result, formulated below.  We use it to discuss \emph{other approaches} to compact quantum groups.

There are several. Among them, we have the work of Dijkhuizen and Koornwinder (\cite{D-K}, 1994). Also Kirchberg presented an approach to compact quantum groups in Oberwolfach (unpublished, 1994). They are essentially based on (or equivalent with) the following property.

\prop\label{prop:3.10}
Assume that $(\mathcal A,\Delta)$ is a Hopf $^*$-algebra with a \emph{positive} integral. The GNS representation of $\mathcal A$ induced by this functional yields bounded operators. The norm closure is a C$^*$-algebra $A$ and $\mathcal A$ can be viewed as a $^*$-subalgebra of $A$. The coproduct $\Delta$ has a unique extension to $A$. The result  is a compact quantum group $(A,\Delta)$ in the sense of Definition \ref{defin:3.6} of Woronowicz.
\eprop

We will not consider these approaches further in these notes. We just learn from them that it is possible to give simple characterizations of compact quantum groups within a purely algebraic framework. This is a consequence of the two results Theorem \ref{thm:3.9} and Proposition \ref{prop:3.10}, formulated above. 

\subsection{Discrete quantum groups}
As we mentioned before, discrete quantum groups appear for the first time in the work of Podle\'s and Woronowicz (\cite{P-W}, 1990). They study discrete quantum groups as duals of compact quantum groups. In particular, properties of discrete quantum groups are derived from properties of compact quantum groups. 

A first independent approach to discrete quantum groups is found in the paper of Effros and Ruan (\cite{E-R}, 1994). Their approach is purely algebraic and the work is in fact very similar to the work of  Dijkhuyzen and Koornwinder  (\cite{D-K}, 1994). From their point of view, Effros and Ruan study discrete quantum groups, while Koornwinder and Dijkhuyzen think of their notion as compact quantum groups. The confusion illustrates Remark \ref{opm:2.6} in Section \ref{s:fqgrps}. In the case of a finite group, we argued why the associated quantum group should be the function algebra. The problem with the two approaches here is precisely that properties of the algebra are formulated \emph{together} with properties of the dual for defining the notion. That is the real origin of the confusion. 

In fact, the treatments of Effros and Ruan on the one hand and of Dijkhuyzen and Koornwinder on the other hand are close to viewing discrete quantum groups as duals of compact quantum groups, as done by Podle\'s and Woronowicz. 
In our work on discrete quantum groups (\cite{VD3n}, 1996), obtained independently around the same time, we approach discrete quantum groups in a more consequent independent way. We first look for a property of a discrete space that we can quantize so as to obtain a notion of a \emph{discrete quantum space}. We define an appropriate notion of a coproduct and finally we impose conditions to make this coproduct group-like.  

Before we give the definition found in \cite{VD3n} and comment on it, we need to recall the notion of  a multiplier Hopf $^*$-algebra as found in \cite{VD1}.

\defin\label{defin:3.11}
Let $\mathcal A$ be a non-degenerate $^*$-algebra over the field of complex numbers. Consider the tensor product $\mathcal A\ot \mathcal A$ and its multiplier algebra $M(\mathcal A\ot \mathcal A)$. A coproduct $\Delta$ on $\mathcal A$ is a $^*$-homomorphism from $\mathcal A$ to $M(\mathcal A\ot \mathcal A)$ satisfying coassociativity. The pair $(\mathcal A,\Delta)$ is a called a multiplier Hopf $^*$-algebra if the linear maps $T_1$ and $T_2$, defined on $\mathcal A\ot \mathcal A$ by
\begin{equation*}
T_1(a\ot b)=\Delta(a)(1\ot b)
\tussenen
T_2(a\ot b)=(a\ot 1)\Delta(b),
\end{equation*}
have range in $\mathcal A\ot \mathcal A$ and are bijective from this space to itself.
\edefin

We refer to \cite{VD1} for the notion of the multiplier algebra of a non-degenerate algebra, as well as for the notion of coassociativity for coproducts with values in the multiplier algebra. 
\snl
Now we recall the definition of a discrete quantum group as given in \cite{VD3n}.

\defin
Let $\mathcal A$ be a direct sum of finite-dimensional full matrix algebras over the field of complex numbers. A coproduct is a non-degenerate $^*$-homomorphism $\Delta$ from $\mathcal A$ to the multiplier algebra $M(\mathcal A\ot \mathcal A)$ of the tensor product  $\mathcal A\ot \mathcal A$ satisfying coassociativity. The pair $(\mathcal A,\Delta)$ is called a \emph{discrete quantum group} if it is a multiplier Hopf $^*$-algebra. 
\edefin

The direct sum here is in an algebraic sense, not a topological one. We add a couple of important remarks about this definition.

\opm
i) Strictly speaking, we should consider the C$^*$-algebraic direct sum $A$ of these matrix algebras, together with a coproduct $\Delta$ from $A$ to the multiplier algebra of the C$^*$-tensor product $A\ot A$ of $A$ with itself. \\
ii) Then we can consider the maps $T_1$ and $T_2$ on the C$^*$-tensor product. They should still be injective and now have dense range. \\
iii) It is easy to see that the (algebraic) direct sum and the C$^*$-algebraic direct sum determine each other. The same is true for the conditions on the canonical maps $T_1$ and $T_2$ are the same. This is however slightly more difficult to show, but on the other hand, intuitively clear.\\
iv) Although not completely in accordance with the main ideas formulated before, we prefer the above definition because we have the notion of a multiplier Hopf algebra available and well understood. 
\eopm	

It is shown in \cite{VD3n} that integrals on discrete quantum groups exist and are unique. This makes it possible to construct the dual. It turns out to be a compact quantum group as discussed before. Moreover, the dual of this dual compact quantum group - now defined as in the work of Podle\'s and Woronowicz - is the original discrete quantum group. 
\snl
We want to emphasize that in \cite{VD3n}, properties of discrete quantum groups are obtained without the knowledge that they are the dual of a compact quantum group. The same is of course true for the construction of the dual compact quantum group of a discrete quantum group. In this sense, the treatment of discrete quantum groups in \cite{VD3n} should be regarded as the right one of discrete quantum groups. We also claim that discrete quantum groups are better understood if considered first in their own right, not referring to the dual compact quantum groups. 
\snl
There are however different opinions about this. 
\snl
The purely algebraic notions of a compact quantum group and a discrete quantum group, as discussed above, fit very well into a more general class, the algebraic quantum groups. They are discussed in a separate section for various reasons we will explain.

\section{\hspace{-17pt}. Algebraic quantum groups }\label{s:alg} % \input artikel4.tex
 
In this note, we use the following terminology. We will motivate our choice.

\defin\label{defin:4.1}
An \emph{algebraic quantum group} is a multiplier Hopf $^*$-algebra with positive integrals.
\edefin

\opm
i) This notion should not be confused with another notion of algebraic quantum groups. In \cite{B-G} the term is used for quantized coordinate rings of algebraic groups. The idea in \cite{B-G} as well as here is that quantum groups are studied within an algebraic framework. This in contrast with the locally compact quantum groups we discuss further where we have operator algebraic concepts.\\
ii) Sometimes any multiplier Hopf algebra with integrals is called an algebraic quantum group. But here it is quite natural to require the \emph{existence of an involution} and \emph{positive} integrals. We will see later (see Proposition \ref{prop:4.9a} below) that this implies that the underlying algebra is an operator algebra. And of course, this is what we need on our way to locally compact quantum groups.
\eopm

\subsection{Properties of algebraic quantum groups}

Any compact quantum group is an algebraic quantum group if we consider the underlying Hopf $^*$-algebra with its integral. We refer to Theorem \ref{thm:3.9} and Proposition \ref{prop:3.10}. In fact, an algebraic quantum group is a compact quantum group if and only the algebra is unital. Similarly any discrete quantum group is an algebraic quantum group if we consider the algebraic direct sum of the building blocks. It has a cointegral. Conversely, an algebraic quantum group with a cointegral is a discrete quantum group. These properties characterize the compact and discrete quantum groups among the more general class of algebraic quantum groups.
\snl
The reader must be aware of the fact that the class of algebraic quantum groups, although it contains compact and discrete quantum groups, and it is self-dual, is not big enough to cover all locally compact quantum groups. If e.g.\ we think of a locally compact group as a locally compact quantum group by associating the pair $(C(G),\Delta)$ as in Definition \ref{defin:3.1}, then it is an algebraic quantum group if and only if $G$ has a compact open subgroup (see \cite{L-VD}).
\snl
Now we recall the main properties of algebraic quantum groups. This will help us to move from algebraic quantum groups to locally compact quantum groups. First we have the following result about integrals.

\prop
If $\varphi$ is a positive left integral, then $\varphi\circ S$ is a right integral and again positive.
\eprop

That $\varphi$ is positive means $\varphi(a^*a)\geq 0$ for all $a\in A$. The positivity of $\varphi\circ S$ is not a trivial result because in general the antipode $S$ is not a $^*$-map. 
Integrals are unique up to a scalar. We will in what follows choose a left integral $\varphi$ and take $\varphi\circ S$ for the right integral $\psi$.

\prop
i) The integrals are faithful.\\
ii) There are automorphisms $\sigma_\varphi$ and $\sigma_\psi$ of $A$ satisfying and characterized by
\begin{equation*}
\varphi(ab)=\varphi(b\sigma_\varphi(a))
\tussenen
\psi(ab)=\psi(b\sigma_\psi(a))
\end{equation*}
for all $a,b\in A$.\\
iii) There is an invertible group-like multiplier $\delta$ in $M(A)$ satisfying and characterized by $\psi(a)=\varphi(a\delta)$ for all $a$.
\eprop

The automorphisms $\sigma_\varphi$ and $\sigma_\psi$ are called \emph{modular automorphisms} and the property is referred to as the KMS-condition (Kubo Martin Schwinger condition) in this context. The terms find their origin in the operator algebra techniques used in mathematical physics. If a finite-dimensional algebra has a faithful linear functional, these automorphisms automatically exist. And algebraists will use the term Nakayama automorphism (see e.g.\ \cite{Na}), in fact for the inverse of what we call the modular automorphism. Observe however that for infinite-dimensional algebras, the existence of a faithful functional is not sufficient to have modular automorphisms. So the above property ii) does not follow from property i).
\snl
The element $\delta$ is called the \emph{modular element}. Here the terminology comes from the theory of locally compact groups. For a locally compact group, the modular function is used to express the right Haar measure in terms of the left Haar measure. This is precisely what appears in item iii) of the above proposition.
\snl
There are many relations between  these objects and the counit and antipode, as well as among these objects themselves. We refer to \cite{VD3}. See also \cite{VD1} for the general theory of multiplier Hopf algebras.

\opm
In the non-involutive case, there is a scalar $\tau$ satisfying $\varphi(S^2(a))=\tau\varphi(a)$ for all $a$. For multiplier Hopf $^*$-algebras with positive integrals, it turns out that this scalar is trivial. The problem was formulated in \cite{Ku-VD} and solved in Theorem 3.4. of \cite{DC-VD1}. It is the only feature where algebraic quantum groups lack the generality of the structure of locally compact quantum groups where this scalar can be non-trivial. We will come back to this later. See e.g.\ Subsection \ref{ss:5.4}, in particular the remarks leading up to Remark \ref{rem:5.21}. 
\eopm

\subsection{Duality for algebraic quantum groups}
The dual of an algebraic quantum group can be defined and turns out to be again an algebraic quantum group. More precisely, we have the following duality.

\prop\label{prop:4.5}
Let $\widehat A$ be the subspace of the linear dual space $A'$ of $A$ of elements of the form $\varphi(\,\cdot\, a)$ where $a\in A$. The coproduct on $A$ yields a product on $\widehat A$ and if we define $\omega^*$ for $\omega\in \widehat A$ by $\omega^*(a)=\overline{\omega(S(a)^*)}$, it is a $^*$-algebra. The adjoint of the product defines a coproduct $\widehat\Delta$ on $\widehat A$ making it into a multiplier Hopf $^*$-algebra. If we define $\widehat\psi(\omega)=\varepsilon(a)$ when $\omega=\varphi(\,\cdot\,a)$ we get a positive right integral on the dual $\widehat A$. Hence $(\widehat A,\widehat \Delta)$ is again an algebraic quantum group.
\snl
The dual of $(\widehat A,\widehat\Delta)$ is canonically identified with the original multiplier Hopf $^*$-algebra. 
\eprop

The duality of algebraic quantum groups includes the duality between compact and discrete quantum groups. If e.g. for a compact quantum group, we consider the underlying Hopf $^*$-algebra as an algebraic quantum group, and if we take the dual in the sense of the above proposition, we arrive at a discrete quantum group. And vice versa. Sometimes, it is easier to deal with the duality between compact and discrete quantum groups when it is considered as a special case of the duality of algebraic quantum groups.

\opm
By defining the dual as in the above proposition, we actually make a choice that is different from what is common in the theory of locally compact quantum groups. The point is that we define the coproduct on the dual by the formula
\begin{equation*}
\widehat \Delta(\omega)(a\ot b)=\omega(ab)
\end{equation*}
for $\omega\in \widehat A$ and $a,b\in A$. This is most natural in Hopf algebra theory (and by extension also for multiplier Hopf algebras) but in the operator algebra approach to quantum groups, rather the flipped version of this coproduct is used. This has many consequences for the concrete forms of the various formulas. One thing is that the antipode on the dual is replaced by its inverse for this choice. Also the dual of the dual is no longer the same as the original one. We have to repeat this procedure and we arrive at the original only after taking the dual four times. 
\eopm

\subsection{The Larson Sweedler theorem. A source of inspiration}\label{ss:LST}
As we will see in Section \ref{s:lcqgrps},  the theory of locally compact quantum groups in the operator algebra setting uses the existence of the Haar weights to develop the theory. This may seem somewhat strange. However, a theorem by Larson and Sweedler \cite{La-Sw}, that appeared already in 1969, states, roughly speaking, that the the existence of integrals implies the existence of the antipode. This is precisely what happens also for locally compact quantum groups. For a locally compact quantum group, the construction of the antipode is found in the paper by Kustermans and Vaes \cite{Ku-V2}. The procedure is explained also in \cite{VD9n}. We will discuss this in Section \ref{s:lcqgrps}because it is an important step in the theory.
\snl
Because of this, it is interesting to see what happens in the algebraic setting. After all, the core of this note is to explain \emph{the road from Hopf algebras to locally compact quantum groups}. It helps to understand the treatment of locally compact quantum groups as given e.g. in \cite{VD9n}.
\snl
We formulate the result for multiplier Hopf algebras (see Theorem 2.3 in \cite{VD7n}). We include a summary of the proof. In Section \ref{s:lcqgrps} we will see how these ideas are used to develop the theory of locally compact quantum groups.

\stel\label{thm:4.8}
Let $A$ be an algebra with a non-degenerate product and assume that $\Delta$ is a full and regular comultiplication on $A$.
If there exists a faithful left integral and a faithful right integral, then $(A,\Delta)$ is a regular multiplier Hopf algebra.
\estel

Recall that a coproduct on $A$ is \emph{full} if its legs are all of $A$. It is called \emph{regular} if not only the the canonical maps $T_1$ and $T_2$ (see Definition \ref{defin:3.11})
have range in $A\ot A$, but the same is true for the canonical maps 
associated with the opposite coproduct $\zeta\Delta$ (where $\zeta$ is the flip map). A multiplier Hopf algebra $(A,\Delta)$ is \emph{regular} if also $(A,\zeta\Delta)$ is a multiplier Hopf algebra . Remember that regularity is automatic when $A$ is a $^*$-algebra and $\Delta$ a $^*$-homomorphism.

\bew  
To prove the result, we have to show that the four canonical maps are bijective. We only give the main ideas, for details we refer to \cite{VD7n}.

Injectivity of the  two canonical maps 
\begin{equation*}
T_2:=a\ot a'\mapsto (a\ot 1)\Delta(a')
\tussenen
T_4:=a\ot a'\mapsto \Delta(a)(a'\ot 1)
\end{equation*}
follows from the existence of a faithful left integral. To show this, multiply e.q.\ in the first case with $\Delta(x)$ from the right and apply a left integral on the second leg. Similarly, the injectivity of the other two maps
\begin{equation*}
T_1:=a\ot a'\mapsto \Delta(a)(1\ot a')
\tussenen
T_3:=a\ot a'\mapsto (1\ot a)\Delta(a')
\end{equation*}
is proven by using the faithful right integral. 
\snl
For the surjectivity the idea is the following. Take elements $a,a',b$ in $A$. Then consider 
\begin{equation*}
(\iota\ot\iota\ot \varphi)(\Delta_{13}(a)\Delta_{23}(a')(1\ot b\ot 1)),
\end{equation*}
where $\varphi$ is the left integral. If we apply the canonical map $T_1$ to this element, we get $c\ot b$ where 
$c=(\iota\ot\varphi)(\Delta(a)(1\ot a'))$. Now we use the fullness of the coproduct and the faithfulness of $\varphi$ to obtain that all elements in $A$ are a linear combination of elements like $c$. This will eventually give the surjectivity of $T_1$. Similarly for the three other maps, where in one case, again the left integral is used and in the remaining two cases, the right integral.
\ebew

 Remark in passing that even for the case of finite-dimensional Hopf algebras, the above result is slightly stronger than the original one in \cite{La-Sw} while the proof is simpler.
\snl
This result is \emph{not sufficient for our purposes}. We know that the antipode exist when the canonical maps are bijective from the general theory of multiplier Hopf algebras. But in this case we want an explicit construction for the antipode. 
\snl
Recall that for a multiplier Hopf algebra, the counit and the antipode are obtained from the inverse of the canonical map $T_1$ with the formulas
\begin{equation*}
mT_1^{-1}(p\ot q)=\varepsilon(p)q
\tussenen
(\varepsilon\ot\iota)T_1^{-1}(p\ot q)=S(p)q
\end{equation*}
for all $p,q\in A$. As before, we use $m$ for the multiplication map (see \cite{VD1}). If we apply this for the result obtained above, we recover the well known formula
\begin{equation*}
S((\iota\ot\varphi)(\Delta(a)(1\ot a')))=(\iota\ot\varphi)((1\ot a)\Delta(a'))
\end{equation*}
for all $a,a'\in A$.
\snl
Now suppose that we try to use this formula to define the antipode. The first problem then is to show that it is well-defined. Let us formulate it under the form of a proposition below.

\prop\label{prop:4.9b}
Assume that we have finitely many  elements $a_i,a'_i$ such that 
\begin{equation*}
\sum_i (\iota\ot\varphi)(\Delta(a_i)(1\ot a'_i))=0.
\end{equation*}
Then 
\begin{equation*}
\sum_i (\iota\ot\varphi)((1\ot a_i)\Delta(a'_i))=0.
\end{equation*}
\eprop

\bew
We have seen above that 
\begin{equation*}
(\iota\ot\varphi)(\Delta(a_i)(1\ot a'_i))\ot 1=T_1((\iota\ot\iota\ot\varphi)\Delta_{13}(a_i)\Delta_{23}(a'_i))
\end{equation*}
for all $i$. Now we use the right integral to get the injectivity of $T_1$ as before. It follows that 
\begin{equation*}
\sum_i (\iota\ot\iota\ot\varphi)((\Delta_{13}a_i)\Delta_{23}(a'_i))=0
\end{equation*}
and if we apply $\varepsilon$ on the first factor we get the desired result.
\ebew

The faithfulness of the right integral is needed to have $S$ well-defined while the faithfulness of the left integral is needed to have it everywhere defined. Similar arguments give that $S$ is also bijective.
\snl
There is still \emph{another characterization of the antipode}. It is based on the following result. Essentially, this appears already in the original paper on multiplier Hopf algebras, see Lemma 5.4 and Lemma 5.5 in \cite{VD1}.

\prop\label{prop:4.10a}
Let $(A,\Delta)$ be a multiplier Hopf algebra. Given an element $a\in A$, for all $a'\in A$ we  write
\begin{equation*}
a\ot a'=\sum_i \Delta(p_i)(1\ot q_i)
\end{equation*}
where $p_i,q_i$ are in $A$. Then we have
\begin{equation*}
(S(a)\ot 1)\Delta(a')=\sum_i (1\ot p_i)\Delta(q_i).
\end{equation*}
\eprop

By the bijectivity of the canonical map $T_1$ we know that we can write $a\ot a'$ as above. Now we use the \emph{Sweedler notation} to write 
\begin{align*}
a\ot a'&=\sum_{(a)} a_{(1)}\ot \varepsilon(a_{(2)})a'\\
&=\sum_{(a)} a_{(1)}\ot a_{(2)}S(a_{(3)})a'\\
&=\sum_{(a)} \Delta(a_{(1)})(1\ot S(a_{(2)})a').
\end{align*}
Then 
\begin{align*}
\sum_{(a)} (1\ot a_{(1)})\Delta(S(a_{(2)})a')
&=\sum_{(a)} S(a_{(3)})\ot a_{(1)}S(a_{(2)})\Delta(a')\\
&=\sum_{(a)} (S(a_{(2)}\ot \varepsilon(a_{(2)}))\Delta(a')=(S(a)\ot 1)\Delta(a').
\end{align*}

One can use this formula to define the antipode. Again the integrals are needed to show that it is well-defined and bijective. We leave the details to the reader and formulate it as an exercise below. This is essentially the method used in \cite{VD9n}. We will say more about this in a subitem on the treatment of locally compact quantum groups in Section \ref{s:lcqgrps}.
\snl
For completeness, we refer also to \cite{Va-VD2} where a similar method is used for C$^*$-Hopf algebras.

\oef
i) Complete the proof of Theorem \ref{thm:4.8}, using the ideas given above. \\
ii) Prove the existence and properties of the antipode, based on the result of Proposition \ref{prop:4.9b}.\\
iii) Finally, use the ideas from Proposition \ref{prop:4.10a} to construct the antipode and again obtain its properties with this method. 
\eoef

Solving this exercise will help to understand what happens with the construction of the antipode for locally compact quantum groups in Section \ref{s:lcqgrps}.

\subsection{The analytic structure of algebraic quantum groups}
With algebraic quantum groups, we get as close as possible to locally compact quantum groups within a purely algebraic context. It is quite remarkable that there is an analytic structure. We recall the following result (Theorem 3.5 from \cite{DC-VD1}).

\stel
Let $(A,\Delta)$ be an algebraic quantum group. Then $A$ is spanned by elements which are simultaneously eigenvectors for $S^2$, $\sigma_\varphi$, $\sigma_\psi$ and left and right multiplication by $\delta$. Moreover, the eigenvalues of these actions are all positive.
\estel

As a consequence, elements in $A$ are all analytic for these actions. We have e.g.\ that there exists a one-parameter group of $^*$-automorphisms $(\sigma_t)_{t\in \mathbb R}$ of $A$ so that all elements of $A$ are analytic and have the property that $\sigma_{-i}=\sigma_\varphi$. Similarly, there is a one-parameter group of $^*$-automorphisms $(\tau_t)_{t\in \mathbb R}$ of $A$, again all elements of $A$ are analytic and $\tau_{-i}=S^2$
\snl
Essentially, these results are already found in \cite{Ku}. In that paper, they are obtained via the passage from algebraic quantum groups to operator algebraic quantum groups. In \cite{DC-VD1}, the results are proven directly, using purely algebraic methods. Along the way, it is shown that the scaling constant $\nu$, satisfying $\varphi\circ S^2=\nu\varphi$ when we have a multiplier Hopf algebra with a left integral $\varphi$, is equal to $1$ if it is an algebraic quantum group in the sense of this note, i.e.\ if we have a multiplier Hopf $^*$-algebra with a positive left integral.
\snl
At this level, we already have a polar decomposition of the antipode. Here it reads as 
$S=R\tau_{-i/2}$ 
where $R$ is a $^*$-anti-automorphism, called the \emph{unitary antipode}. Remark that again, implicitly, a choice is made. It complies with the notions used for locally compact quantum groups in \cite{Ku-V2}. In the work of Masuda, Nakagami and Woronowicz \cite{M-N,M-N-W} however a different choice for the scaling automorphism group is made and consequently the polar decomposition of the antipode appears as 
$S=R\tau_{i/2}$. 

\subsection{From algebraic quantum groups to operator algebraic quantum groups}
Consider an algebraic quantum group $(A,\Delta)$ and let $\psi$ be a right integral. The GNS-space associated with $\psi$ is a Hilbert space, denoted by $\mathcal H_\psi$, together with an injective linear map $\Lambda:A\to \mathcal H_\psi$ with dense range and such that $\langle\Lambda(a),\Lambda(b)\rangle=\psi(b^*a)$ for all $a,b$. For the GNS-representation of $A$ we find the following property.

\prop\label{prop:4.9a}
There is a non-degenerate $^*$-representation  $\pi_\psi$ of $A$ by \emph{bounded} operators on $\mathcal H_\psi$ satisfying $\pi_\psi(a)\Lambda(b)=\Lambda(ab)$ for all $a,b$.
\eprop

Most of this is obvious, except for the fact that the map $\Lambda(b)\mapsto \Lambda(ab)$ is a bounded map. The proof  is not difficult but needs more than just the property of a positive linear functional on the $^*$-algebra $A$. It is indeed not difficult to find an example of a $^*$-algebra with a positive faithful linear functional on it so that the GNS representation is not by bounded operators.
\snl
The proof needs the invariance of $\psi$ and it is related with the following result. To make the formulation somewhat easier,  we let $A$ act directly on the $\mathcal H_\psi$ and so we write $a\xi$ for $\pi_\psi(a)\xi$ when $\xi\in\mathcal H_\psi$ and $a\in A$. 

\prop\label{prop:4.9}
 There is a unitary operator $V$ on $\mathcal H_\psi\ot\mathcal H_\psi$ satisfying
\begin{equation*}
V(\Lambda(a)\ot \Lambda(b)=(\Lambda\ot\Lambda)(\Delta(a)(1\ot b))
\end{equation*}
for all $a,b\in A$.
\eprop

Again the proof is easy. By definition the canonical map $a\ot b\mapsto \Delta(a)(1\ot b)$ is a bijection from $A\ot A$ to itself. So the map $V$ is well-defined on the dense subspace $\Lambda(A)\ot\Lambda(A)$ of $\mathcal H_\psi\ot \mathcal H_\psi$. It is isometric because $\psi$ is right invariant. Finally, because the canonical map is a bijection, the isometric operator $V$ will be a unitary. 
\snl
Remark that we have a similar result for general locally compact quantum groups. In that case, it is also easy to get that this operator is isometric, but not that it is a unitary. We will come back to this later. See Problem \ref{probl:6.10a}.
\snl
The unitary $V$ is a multiplicative unitary and it is the basic ingredient that allows us to pass from the algebraic quantum group to the associated locally compact quantum group. This is done in \cite{Ku-VD} but a more recent  approach to this procedure is found in a forthcoming paper \cite{VD10n}. 

\opm\label{opm:4.10}
i) It is also possible to use the left integral to produce a multiplicative unitary. For this one considers the GNS space $\mathcal H_\varphi$ and the map $\Lambda_\varphi:A\to \mathcal H_\varphi$ associated with a left integral $\varphi$. Then a unitary $W$ is defined by its adjoint with the formula
\begin{equation*}
W^*(\Lambda_\varphi(a)\ot \Lambda_\varphi(b)=(\Lambda_\varphi\ot\Lambda_\varphi)(\Delta(b)(a\ot 1))
\end{equation*}
for $a,b\in A$.\\
ii) As the reader can notice, this definition is somewhat less natural. First there is the map $a\ot b\mapsto \Delta(b)(a\ot 1)$ that is different but derived from the basic canonical maps used to define a multiplier Hopf algebras, see Definition \ref{defin:3.11}. Secondly, the unitary is  defined in terms of the inverse of this map. \\
iii) These two different conventions are used in the theory. In their work on multiplicative unitaries \cite{Ba-Sk}, Baaj and Skandalis prefer to use the one associated with the right Haar measure while in most other papers on locally compact quantum groups, the left regular representation with $W$ as in i) is used.
\eopm

This is again somewhat annoying, certainly because  multiplicative unitaries play an important role as an intermediate
Indeed, before we start with the discussion of the general locallly compact quantum groups in Section \ref{s:lcqgrps}, we include a larger section on multiplicative unitaries. We discuss various aspects of these.
  
%%%%%% Multiplicative unitaries
  
\section{\hspace{-17pt}. Multiplicative unitaries} \label{s:mu} % \input artikel5.tex

As we saw in Proposition \ref{prop:4.9} at the end of the previous section, multiplicative unitaries arise in a natural way when  moving  from the purely algebraic setting to the operator algebraic approach to quantum groups. Indeed, they play an important role in the theory of locally compact quantum groups as we will see in the next section. Therefore, it is of great importance to fully understand how they pop up in the theory. The best way to do this, is to look first at the case of a finite-dimensional Hopf $^*$-algebra.  Further it is interesting to consider algebraic quantum groups because still for them it is possible to treat multiplicative unitaries in a purely algebraic manner.

\subsection{Multiplicative unitaries in the algebraic setting}
We begin with the case of a finite-dimensional Hopf $^*$-algebra $A$. Denote the dual by $B$ and use the pairing notation. So we write $\langle a,b\rangle$ for the value of a linear functional $b$  on an element $a$ in $A$.  The tensor product $B\ot A$ is again a Hopf $^*$-algebra. Its dual can be identified with $A\ot B$. We consider the natural associated pairing of $B\ot A$ with $A\ot B$ given by
$$\langle b\ot a,a'\ot b'\rangle=\langle a',b\rangle\langle a,b'\rangle$$
 for $a,a'\in A$ and $b,b'\in B$.
\snl
Because all spaces here are finite-dimensional, we can define  $V$ in $B\ot A$ by 
\begin{equation}
\langle V,a\ot b\rangle=\langle a,b\rangle\label{eqn:5.1a}
\end{equation}
where $a\in A$ and $b\in B$. We have the following properties.

\prop\label{prop:5.1}
The element $V$ is a unitary in the $^*$-algebra $B\ot A$. We have
\begin{equation}
(\varepsilon\ot\iota)V=1 
\tussenen
(\iota\ot\varepsilon)V=1 \label{eqn:5.1}
\end{equation}
when we use $\varepsilon$ for the counit on $A$ as well as for the one on the dual $B$. For the antipode we find
\begin{equation}
(S\ot\iota)V=V^* \tussenen (\iota\ot S)V=V^* \label{eqn:5.2}
\end{equation}
where again we use $S$ both for the antipode on $A$ and on $B$. Finally, using $\Delta$ for the coproduct on $A$ as well as on $B$ we find
\begin{equation}
(\Delta\ot\iota)V=V_{13}V_{23}
\tussenen
(\iota\ot\Delta)V=V_{12}V_{13}. \label{eqn:5.3}
\end{equation}
\eprop

We use the \emph{leg-numbering notation}. We have $V_{23}=1\ot V$ in $B\ot B\ot A$ and $V_{12}=V\ot 1$ in $B\ot A\ot A$. The element $V_{13}$ belongs to $B\ot B\ot A$ in the first equation and to $B\ot A\ot A$ in the second equation. It is obtained by inserting $1$ in the middle factor, the identity in $B$ in the first case and in $A$ for the second formula.
\snl
The element $V$ is thought of as \emph{the duality}. The formulas in Proposition \ref{prop:5.1} are nothing more than reformulations of the properties of a pairing between Hopf $^*$-algebras. 

\oef 
Verify that the properties of a pairing of Hopf $^*$algebras are reformulated in terms of the duality $V$ as in Proposition \ref{prop:5.1}. To show that $V$ is a unitary one uses the basic properties of the antipode.
\eoef

There is a \emph{similar result} for algebraic quantum groups. In this case $V$ is defined in the multiplier algebra $M(B\ot A)$ and the formulas stated in the proposition have to be formulated in multiplier algebras. 

\oef
i) Consider an algebraic quantum group $A$ and denote also here its dual $\widehat A$ by $B$. 
 First one defines $V$ as an element in the dual of $A\ot B$ as above in Equation (\ref{eqn:5.1a}). To show that it actually is defined as an element in the multiplier algebra $M(B\ot A)$, one first has to extend the pairing. This has been done in various papers on multiplier Hopf algebras and algebraic quantum groups. See e.g.\ Section 3 in \cite{De-VD} where this is done in detail for algebraic quantum hypergroups. 
 \\
ii) This is in fact the main issue. Once this is done, proving the formulas is similar as for the case of finite dimensional Hopf $^*$-algebras.\\
iii) A new treatment will be found  in \cite{VD10n}.
\eoef

In order to view $V$ as a \emph{multiplicative unitary}, we need the following result. We are again in the finite-dimensional case and once more, we use  the Sweedler notation in what follows.

\prop\label{prop:5.2}
Let $A$ act on itself from the left by multiplication. We also let $B$ act from the left on $A$ by the formula
\begin{equation*}
(b,a)\mapsto b\tr a:= \sum_{(a)}\langle a_{(2)},b\rangle a_{(1)}.
\end{equation*}
Then $B\ot A$ acts from the left on $A\ot A$. Using  $\tr$ for all these different actions, we find 
\begin{equation*}
V\tr(a\ot a')=\Delta(a)(1\ot a') 
\end{equation*}
for all $a,a'\in A$.
\eprop

The proof is simple. If we write $V=\sum_i v_i\ot w_i$ with $v_i\in B$ and $w_i\in A$ we find
\begin{equation}
V\tr (a\ot a')
=\sum_i v_i\tr a\ot w_ia' 
=\sum_{i,(a)} \langle a_{(2)},v_i\rangle a_{(1)}\ot w_ia'.\label{eqn:5.5a}
\end{equation}
Because $V$ is the duality, we have
$\sum_{i} \langle p,v_i\rangle w_i=p$ for all $p\in A$ and hence the last expression in Equation (\ref{eqn:5.5a}) is equal to $\sum_{(a)} a_{(1)}\ot a_{(2)}a'$. This proves the result.
\snl
Now we have $V(a\ot 1)=\Delta(a)V$ and $V(a\ot 1)V^*=\Delta(a)$ for all $a\in A$ as \emph{equalities of linear operators} on $A\ot A$. These formulas however have also an \emph{interpretation using the Heisenberg algebra}. This is based on the following result.

\prop\label{prop:5.3}
For $a,a'\in A$ and $b\in B$ we have
\begin{equation*}
ba=\sum_{(a),(b)}\langle a_{(2)} ,b_{(1)} \rangle a_{(1)}b_{2}.
\end{equation*}
Also this formula has to be viewed as an equation of linear maps on $A$. 
\eprop
Indeed, if also $a'\in A$ we find
\begin{align*}
b\tr (aa')
&=\sum_{(a),(a')}\langle a_{(2)}a'_{(2)},b \rangle a_{(1)}a'_{(1)}\\
&=\sum_{(a),(a'),(b)}\langle a_{(2)},b_{(1)} \rangle\langle a'_{(2)},b_{(2)} \rangle a_{(1)}a'_{(1)}\\
&=\sum_{(a),(b)}\langle a_{(2)},b_{(1)} \rangle a_{(1)}(b_{(2)}\tr a').
\end{align*}

This leads us to the following definition.

\defin\label{defin:5.6a}
The formula in the previous proposition is called the \emph{Heisenberg commutation relation}. 
The \emph{Heisenberg algebra} is the algebra generated by $A$ and $B$ subject to the Heisenberg commutation relations.
\edefin

The terminology is justified because in the case of the group $\mathbb R$ and the locally compact quantum group associated to it, the  commutation rules in Proposition \ref{prop:5.3} will essentially be the same as the Heisenberg commutation rules (or its Weyl equivalent forms) if interpreted in the right way. 
\snl
We denote the Heisenberg algebra simply by $AB$. It has a faithful left action on $A$ induced by the actions of $A$ and of $B$. That we have a left action is obvious. To show that it is faithful requires an argument.

\opm
One should be careful with this terminology. Indeed, the Heisenberg algebra in itself is trivial. For a finite-dimensional Hopf algebra it is a full matrix algebra. And in the more general case of an algebraic quantum group, it is an infinite matrix algebra. It only depends on the linear spaces $A$ and $B$ and their pairing, but not on the algebra structures. See e.g.\ Proposition 6.6 in \cite{D-VD-Z}. For this reason, it is better to speak about the Heisenberg algebra $AB$ as the algebra \emph{together} with the embeddings of $A$ and of $B$ in (the multiplier algebra of) $AB$.  
\eopm

Now the formula $V(a\ot 1)V^*=\Delta(a)$ is valid in the algebra $AB\ot A$. And if we combine this formula with the formulas involving the coproduct in Proposition \ref{prop:5.1} we naturally arrive at the Pentagon equation
\begin{equation}
V_{12}V_{13}V_{23}=V_{23}V_{12}, \label{eqn:5.4}
\end{equation}
a formula that holds here in the algebra $B\ot AB\ot A$.
\snl
All the above results are still valid for the duality of algebraic quantum groups $(A,\Delta)$. As we mentioned already, in that case the element $V$ is  an element in the multiplier algebra $M(B\ot A)$ where now $B$ is the dual $\widehat A$ of $A$ as in Proposition \ref{prop:4.5}. Also the Heisenberg commutation relations can be formulated and adapted forms of the results in Propositions \ref{prop:5.2} and \ref{prop:5.3} are true. The Pentagon equation is now valid in the multiplier algebra of $B\ot AB\ot A$. The results are present in \cite{Ku-VD} but a more direct and simpler approach will appear in the forthcoming paper \cite{VD10n}.

\opm
i) If we take the other convention for defining the coproduct on the dual (by flipping it),  some of the formulas formulated in Proposition \ref{prop:5.1} will change. The formulas (\ref{eqn:5.1}) involving the counit will remain the same but the second formula of (\ref{eqn:5.2}) will change to $(S^{-1}\ot\iota)V=V^*$ because the antipode on $B$ is replaced by its inverse. This means e.g.\ that then we have $(S\ot S)V=V$. Also the first formula of (\ref{eqn:5.3}) will be replaced by $(\Delta\ot\iota)V=V_{23}V_{13}$. The Pentagon equation (\ref{eqn:5.4}) remains unchanged. This is natural because the definition of $V$ does not depend on the definition of the coproduct on the dual. Remark however that the Heisenberg commutation relations do change.\\
ii) One of the advantages of the operator algebraic convention over the algebraic one is the formula $(S\ot S)V=V$. Also the multiplicative unitary for the dual has a nice relation with the original.
 \\
iii) But still, this is a rather annoying difference if you approach locally compact quantum groups from the Hopf algebra theory.
\eopm

\oef
Most of the previous statements are easy to verify in the case of a finite-dimensional Hopf $^*$-algebra. On the other hand, with a suitable reformulation, they are also true for algebraic quantum groups. The arguments are essentially the same, but some care is needed. Indeed, generally speaking, if results are true for the Hopf algebras and if the formulations make sense for multiplier Hopf algebras, they tend to also be true for those. It only takes more careful work to translate the proofs. The above is an excellent set of results to practice this. 
\eoef

We have seen how the canonical map $a\ot a'\mapsto \Delta(a)(1\ot a')$ is related with the duality $V$ and that the Pentagon equation for $V$ has an interpretation in the algebra $B\ot AB\ot A$ (or its multiplier algebra in the case of algebraic quantum groups). But of course, the Pentagon equation satisfied by $V$ can also be seen as a result of linear maps on $A\ot A\ot A$ with $V$ now defined as the map  $a\ot a'\mapsto \Delta(a)(1\ot a')$ from $A\ot A$ to itself. 
\snl
However, in order to call it a unitary in that setting, we need the right integral and its GNS space to lift $V$ to a unitary operator as in Proposition \ref{prop:4.9}. This is the more common approach to multiplicative unitaries. We see this in the next subsection.

\subsection{The origin of multiplicative unitaries}
The first appearance of what is now called a multiplicative unitary is  the unitary operator $V$ on $L^2(G\times G)$ for a locally compact group $G$. The Hilbert space is constructed using the right Haar measure on $G$ and the linear map is given by the formula
\begin{equation*}
(Vf)(p,q)=f(pq,q)
\end{equation*}
where $f\in L^2(G\times G)$ and $p,q\in G$.  This is indeed the analogue of the map $a\ot a'\mapsto \Delta(a)(1\ot a')$ we considered before for (multiplier) Hopf algebras, now in terms of the coproduct $\Delta$ on $C(G)$ as given in Definition \ref{defin:3.1} in the case of a compact group. In fact, the unitary $W$ on the Hilbert space $L^2(G\times G)$ constructed with the left Haar measure and defined by
\begin{equation*}
(Wf)(p,q)=f(p,p\inv q)
\end{equation*}
for $f\in L^2(G\times G)$ and $p,q\in G$ is a more familiar object. We recognize here the inverse of the map $a'\ot a\mapsto \Delta(a)(a'\ot 1)$ defined for (multiplier) Hopf algebras. See Proposition \ref{prop:4.9} and Remark \ref{opm:4.10}.

\oef
Verify these statements.
\eoef

\opm
i) The unitary $V$ is thought of as the \emph{right regular representation} whereas $W$ is the \emph{left regular representation} from this point of view.\\
ii) In the literature, it is more common to work with the left Haar measure on a group and consequently to use the unitary $W$ as the multiplicative unitary for further developing the results. On the other hand, we see that from the point of view of Hopf algebra theory, it seems more natural to consider $V$. \\
iii) For this reason, some people will prefer to work with $V$ instead. This is e.g.\ the case in \cite{Ba-Sk}, see Item 1.2.2 in that paper. \\
iv) In the approach to locally compact quantum groups in the work of Kustermans and Vaes \cite{Ku-V1,Ku-V2}, in the beginning of the development, the two are important. See also the more recent treatment in \cite{VD9n}.
\eopm

The unitary operator $W$ has played an important role from the very beginning in the search for a generalization of Pontryagin's duality for non abelian locally compact groups. In that context, it is sometimes called the Kac-Takesaki operator. It also plays an important role in the theory of Kac algebras, developed later. 
\snl
Already as early as in 1979 an attempt was made to develop a theory starting from a multiplicative unitary with some extra properties, see \cite{VH}. The extra properties were inspired by the results obtained for Kac algebras in the work of Enock and Schwartz \cite{E-S1,E-S2}. Remark that at that time, the term multiplicative unitary was not yet introduced.
\snl
This only happened in 1993 with the work of Baaj and Skandalis \cite{Ba-Sk}. In fact, their work was available as a preprint (in different versions) several  years before it actually was published. The theory of multiplicative unitaries in \cite{Ba-Sk} is stronger than the earlier one in \cite{VH} and took advantage of the new developments at that time, with the theory of compact quantum groups by Woronowicz. 
\snl
Recently, another attempt was made to use multiplicative unitaries as a bases for developing the theory (\cite{M-VD2}, 2002). In a way, it is an update of the earlier work of Vanheeswijk (\cite{VH}, 1979). This paper is only available on the arxiv. It was never published, but it does contain some interesting results. We will say more about it later in this section.  

\subsection{Regular multiplicative unitaries}
Recall the definition of a multiplicative unitary on Hilbert spaces (Definition 1.1 in \cite{Ba-Sk}).

\defin\label{defin:5.7}
 Let $\mathcal H$ be a Hilbert space. Denote by $\mathcal H\ot\mathcal H$ the Hilbert space tensor product of $\mathcal H$ with itself. A unitary operator on $\mathcal H\ot\mathcal H$ is called a \emph{multiplicative unitary} if it satisfies the \emph{Pentagon equation}
\begin{equation*}
V_{12}V_{13}V_{23}=V_{23}V_{12}.
\end{equation*}
\edefin
We are now using the leg-numbering notation again. All the factors in this equation are unitaries on the three-fold tensor product $\mathcal H\ot\mathcal H\ot\mathcal H$.
\snl
By the \emph{left leg} of $V$ we mean the subspace of bounded operators on $\mathcal H$ of the form $(\iota\ot\omega)V$ where $\omega$ is an element in $\mathcal B(\mathcal H)_*$, the space of normal linear functionals on  $\mathcal B(\mathcal H)$. The {\emph{right leg} is the space of operators of the form $(\omega\ot \iota)V$ where again $\omega\in\mathcal B(\mathcal H)_*$. It is an immediate consequence of the Pentagon equation that these spaces are subalgebras of $\mathcal B(\mathcal H)$. 
\snl
For reasons that will become clear later, we introduce the following notations.

\notat
We use  $A_0$ and $\widehat A_0$ for the right and the left leg of $V$ respectively and  we denote their  norm closures by $A$ and $\widehat A$.  For the closures w.r.t.\  the $\sigma$-weak operator topology, we write $M$ and $\widehat M$.
\enotat

The spaces $A$ and $\widehat A$ are normed algebras while the spaces $M$ and $\widehat M$ are $\sigma$-weakly closed subalgebras of $\mathcal B(\mathcal H)$. 
\snl
For a multiplicative unitary associated with a locally compact quantum group, these algebras are the group algebras and the function algebras. They are obviously self-adjoint. In general however, they are \emph{not} self-adjoint. It is one of the main issues in the theory of multiplicative unitaries to formulate conditions so that these spaces \emph{are} self-adjoint.
\snl
In the paper of Baaj and Skandalis, there is a result of this type based on the notion of \emph{regularity of the multiplicative unitary}. Regularity of $V$  is formulated in terms of the algebra $C(V)$, defined as the set of operators on $\mathcal H$ of the form $(\iota\ot\omega)(\zeta V)$ where we use $\zeta$ here for the flip map on the tensor product $\mathcal H\ot\mathcal H$ and with again $\omega$ any element of $\mathcal B(\mathcal H)_*$. We refer to Definition 3.3 in \cite{Ba-Sk}. In  Proposition 3.5  of \cite{Ba-Sk} it is shown that for a regular multiplicative unitary, the algebras in Definition \ref{defin:5.7} are all self-adjoint. So $A$ and $\widehat A$ are C$^*$-algebras and $M$ and $\widehat M$ are von Neumann algebras. It is also shown in \cite{Ba-Sk} that $V$ belongs to the multiplier algebra $M(\widehat A\ot A)$ of the minimal C$^*$-tensor product $\widehat A\ot A$ (see Proposition 3.6 in \cite{Ba-Sk}). This implies that $V$ also belongs to the von Neumann algebra tensor product $\widehat M\ot M$, but that result is easier to obtain. 
\snl
It is interesting to see the relation with the Heisenberg algebra as we discussed it in the first subsection of this section. We formulate it as an example.

\voorb i) Let $A$ be a finite-dimensional Hopf algebra. Consider the canonical map $V$ on $A\ot A$ given as before by $V(a\ot a')=\Delta(a)(1\ot a')$. If we also denote the flip map on $A\ot A$ by $\zeta$ we will find
\begin{equation*}
(\zeta V)(a\ot a')=\zeta\sum_{(a)} a_{(1)}\ot a_{(2)}a'=\sum_{(a)} a_{(2)}a'\ot a_{(1)}
\end{equation*}
and if we pair this with an element $b$ of the dual in the second factor we get
\begin{equation*}
((\iota\ot\omega)(\zeta V))a=\sum_{(a)} a_{(2)}a'\langle a_{(1)},b\rangle
\end{equation*}
where $\omega=\langle \,\cdot\,a',b\rangle$. We can write the right hand side as $(a\tl b)a'$ where now
\begin{equation*}
 a\tl b:=\sum_{(a)} a_{(2)}\langle a_{(1)},b\rangle.
\end{equation*}
 We immediately see that this is similar to the left action of the Heisenberg algebra we had in Proposition \ref{prop:5.2}. In fact, the antipode will convert this action precisely to what we had in that proposition. \\
ii) We know that in the case of a finite-dimensional Hopf $^*$-algebra the Heisenberg algebra is the algebra all linear maps on $A$. And in the case of an algebraic quantum group, it is an algebra of finite-rank operators build with the spaces from their pairing. See a remark after Definition \ref{defin:5.6a} earlier in this section.\\
iii) We see from this that the multiplicative unitary associated with an algebraic quantum group as in Proposition \ref{prop:4.9} will automatically be regular.
\evoorb

Remark once more that various choices in the theory of multiplicative unitaries and in particular in their relation with quantum groups are possible. We see this also here, in the above example. This again is another annoyance due to the different approaches.
\snl
After the appearance of \cite{Ba-Sk}, it was shown in \cite{B-S-V} that there are locally compact quantum groups with an associated multiplicative unitary that is not regular. This implied the need for a notion beyond regularity, still sufficient to show that the associated algebras are self-adjoint. 
\snl
Also the antipode gets little attention in the paper of Baaj and Skandalis. Observe that at the time they studied multiplicative unitaries, Kac algebras were known already but locally compact quantum groups as defined more recently were not. 
\snl
There is another property of multiplicative unitaries that guarantees that the algebras $A$ and $\widehat A$ are C$^*$-algebras. It is found in  the work of Woronowicz (\cite{W4}, (1996)) and So\l tan (\cite{S-W}, 2007) on \emph{manageable multiplicative unitaries}. We discuss this in the next subsection. 

\subsection{Multiplicative unitaries and the antipode}\label{ss:5.4}
From Proposition \ref{prop:5.1} in the algebraic context, it is clear what the relation is between the antipode and the multplicative unitary. The antipode on the left leg is determined by $(S\ot\iota)V=V^*$ while the antipode on the right leg by $(\iota\ot S)V=V^*$. 
Hence we see that the  antipode $S$ on e.g.\ the right leg of $V$ will exist if and only if  for any element $\omega\in \mathcal B(\mathcal H)_*$ it is true that  $( \omega\ot\iota)V^*=0$ as soon as $(\omega\ot\iota)V=0$. This condition is necessary and allows us to define $S$ on the right leg of $V$ by
\begin{equation*}
S((\omega\ot\iota)V)=(\omega\ot\iota)V^*
\end{equation*}
for all $\omega\in \mathcal B(\mathcal H)_*$. Similarly for the antipode on the left leg. 
\snl
This however is of \emph{little use}. But if we impose the requirement that not only these maps exist, but also are well-behaved (in a sense we formulate in the definition below), we get a useful property.

\defin\label{defin:5.10} Let $M$ be a von Neumann algebra. 
Assume that $R$ is a linear involutive $^*$-anti-automorphism of $M$ and that $(\tau_t)_{t\in \mathbb R}$ is a $\sigma$-weak operator continuous one-parameter group of $^*$-automorphisms of $M$. We define analytic extensions in an appropriate way. We are particularly interested in the linear operator $\tau_{-i/2}$ on $M$. In general, it is not everywhere defined, but a closed densely defined operator. We now look at the composition $R\circ \tau_{-i/2}$. It will again be a closed densely defined linear map which we will denote by $S$. We say that $S$ has a \emph{polar decomposition} if it arises this way.
\edefin

We will require that elements of the right leg of $V$ belong to the domain $\mathcal D(S)$ of $S$ and that $S((\omega\ot\iota)V)=(\omega\ot\iota)V^*$ for all $\omega\in \mathcal B(\mathcal H)_*$. This is now what we write as $(\iota\ot  S)V=V^*$.
\snl
For an antipode $S$ satisfying this property, it will automatically follow that $S(x)^*$ belongs to the domain of $S$ for all $x$ in that domain and that $S(S(x)^*)^*=x$. We know form Hopf $^*$-algebra theory that this is a natural condition for the antipode. When $S$ has a polar decomposition as in Definition \ref{defin:5.10} above, we will have this property of the antipode if and only if $R$ and $\tau_t$ commute for all $t$. Therefore, it is natural to add this condition above.
\snl
For convenience, we temporarily introduce the following notion.

\defin\label{defin:5.11}
 Assume that $V$ is a multiplicative unitary on $\mathcal H\ot\mathcal H$. Let $\widehat S$ and $S$ be closed linear maps  on $\mathcal B(\mathcal H)$ that admit a polar decomposition as in Definition \ref{defin:5.10}. Assume that the right leg of $V$ belongs to the domain of $S$ and that $(\iota\ot S)V=V^*$. Similarly assume that the left leg of $V$ belongs to the domain of $\widehat S$ and that $(\widehat S\ot \iota)V=V^*$.  Then we say that $V$ is a \emph{multiplicative unitary with manageable  antipodes}.
\edefin

It is proven by So\l tan and Woronowicz   in \cite{S-W} that a \emph{modular multiplicative unitary}, as in Definition  2.1 of \cite{S-W}, is a multiplicative unitary with manageable antipodes as in the above definition, see Theorem 2.2 again in \cite{S-W}. This result generalizes the earlier result by Woronowicz in \cite{W4} for \emph{manageable multiplicative unitaries} as defined in Definition 1.2 of \cite{W4}. 
\snl
Roughly speaking, the property of  modularity for  multiplicative unitary as given in Definition 1.2 of \cite{S-W} is very close to the requirement that it has manageable antipodes as in our definition above \ref{defin:5.11}.
\snl
It is also not so difficult to prove now the following result.

\prop\label{prop:5.12a}
Assume that $V$ is a multiplicative unitary with manageable antipodes.
Then the norm and $\sigma$-weak  operator closures of the legs are self-adjoint.
\eprop

\bew
The idea to show this for the left leg is to use linear functionals $\omega$ in $\mathcal B(\mathcal H)_*$ that are analytic for the one-parameter group $(\tau_t)$ that appears in the polar decomposition of $S$. For such a linear functional $\omega$, we have another element $\omega_1$ in  $\mathcal B(\mathcal H)_*$ so that  $\omega_1(x)=\omega(\tau_{i/2}(x))$ for $x$ in the domain of $\tau_{i/2}$. Then using  $(\iota\ot S)V=V^*$ a straightforward calculation will give that
\begin{equation*}
((\iota\ot\omega)V)^*=(\iota\ot\omega_2)V
\end{equation*}
where $\omega_2$ is defined by $\omega_2(x)=\overline{\omega(R(x)^*)}$ for all $x\in \mathcal B(\mathcal H)$. We use the $^*$-anti-automor\-phism $R$ from the polar decomposition of $S$.
\ebew

The result is found in Theorem 2.2 of  \cite{S-W}. 

\opm
i) As you can see, there is a fundamental difference in the two approaches. The condition of regularity in the work of Baaj and Skandalis \cite{Ba-Sk} is of a different nature than the requirement of having manageable antipodes as in the work of So\l tan and Woronowicz \cite{W4,S-W} described above.\\
ii) The example in \cite{B-S-V} shows that regularity is too restrictive for dealing with locally compact quantum groups. One needs the other approach for that.
\eopm

There is still another approach that is worthwhile mentioning. It builds further on the ideas of the earlier paper \cite{VH}, but now taken into account the recent progress and the new knowledge. It is found in the non-published work, available on the Arxiv (\cite{M-VD2}, 2002). 
\snl
The notion defined in that paper is the following. The condition in \cite{M-VD2} is formulated in terms of the right regular representation $W$ but we will translate it to the framework we are using here.

\defin 
Let $V$ be a multiplicative unitary on $\mathcal H \ot\mathcal H$ and let $A_0$ and $\widehat A_0$ be the right and left leg of $V$ as in Definition \ref{defin:5.7}. Then $V$ is called \emph{trim} if the set $A_0(\widehat A_0)^*$ is $\sigma$-weakly dense in $\mathcal B(\mathcal H)$.
\edefin 

Remark that one can show that the $\sigma$-weak closure of $A_0(\widehat A_0)^*$ is always a subalgebra of $\mathcal B(\mathcal H)$. See Proposition 3.1 and Definition 3.2 in \cite{M-VD2}. The main definition is now Definition 3.4 in \cite{M-VD2}. With the notations used here, it reads as follows.

\defin
Let $V$ be a multiplicative unitary and assume that it is trim as in the above definition. Assume moreover that there are  involutive, self-adjoint  anti-linear operators $J$ and $\widehat J$ on $\mathcal H$ such that $V^*=(J\ot \widehat J)V(J\ot \widehat J)$ and
\begin{equation*}
JMJ\subseteq M'
\tussenen
\widehat J \widehat M\widehat J\subseteq \widehat M'. 
\end{equation*}
Here we use $M$ and $\widehat M$ for the von Neumann algebras \emph{generated by} the legs $A_0$ and $\widehat A_0$ respectively.
\snl
Then the triple $(V,J,\widehat J)$ is called a \emph{quantum frame} on the Hilbert space $\mathcal H$.
\edefin

Also for quantum frames, one can show that the $\sigma$-weak closures of the spaces $A_0$ and $\widehat A_0$ are self-adjoint (see Theorem 3.10 in \cite{M-VD2}. Then the notation for $M$ and $\widehat M$ used in the above definition are compatible with the one used earlier here in Definition \ref{defin:5.7}.
\snl
For a multiplicative unitary of a quantum frame, the antipodes themselves are not obtained, but the unitary antipodes $R$ and $\widehat R$ do exist. They are implemented by the involutions $J$ and $\widehat J$. For a multiplicative unitary with manageable antipodes, we do have the unitary antipode $R$ and $\widehat R$ by definition, but they are not necessarily implemented by elements $J$ and $\widehat J$. This is the extra condition we then have for quantum frames. 
\snl
In Proposition 3.13 of \cite{M-VD2} we find the relation with the Heisenberg algebra. While in Proposition 3.14 of \cite{M-VD2} it is shown that there is a complex number $\lambda$ with modulus $1$ satisfying $J\widehat J=\lambda \widehat J J$. 
\snl
Also observe that any locally compact quantum group gives rise to a quantum frame (see Proposition 3.16 in \cite{M-VD2}). The factor $\lambda$ is related with the scaling constant, see remark 3.17 in \cite{M-VD2}. It is somewhat remarkable that the scaling constant already appears in the theory of quantum frames even though there are no Haar weights around.
\snl
In Section 4 of \cite{M-VD2}, the whole theory is well-illustrated with one of the highly non-trivial examples obtained by Woronowicz, namely the quantum $az+b$-group \cite{W5}. This is an example where the scaling constant is non-trivial.

\opm\label{rem:5.21}
The quantum $az+b$-group is one of the many examples of this kind obtained by Woronowicz. The general procedure is to construct the multiplicative unitary and prove its properties. In order to fit these examples in the theory of locally compact quantum groups however, one needs to have the Haar weights. Although there is no general theory where the Haar weights are constructed from the axioms, in most of these examples it is fairly obvious what the Haar weights have to be. It is also not a problem to show that they are appropriate weights. The main difficulty lies in proving invariance. The reason why this is not easy has to do with the fact that the algebras are defined by a dense subalgebra and that invariance of the weight on this dense subalgebra is not sufficient to have over all invariance. Recall that weights can be different even when they agree on a dense subalgebra. For some of these examples, the Haar weights are found in \cite{VD6}. In that paper, a general technique is described for proving invariance in such cases. We come back to this in the next subsection.
\eopm

\subsection{From multiplicative unitaries to locally compact quantum groups}
The papers by Woronowicz and So\l tan \cite{W4,S-W} bear the title \emph{From multiplicative unitaries to quantum groups}. This is somewhat misleading. It all depends on what you consider to be a quantum group. And as we have seen already, the notion is used for many different things. Since in this paper, we are aiming at locally compact quantum groups, we insist on having the analogues of the Haar measures. Therefore we now address the problem of finding the Haar weights from the properties of a multiplicative unitary.
\snl

In \cite{VD6} we have used a procedure to obtain the Haar weights for the quantum groups associated with some of the multiplicative unitaries as obtained by Woronowicz. The underlying idea is very simple and we explain it below in the finite-dimensional setting. 

\prop
Take a finite-dimensional Hopf $^*$-algebra $A$, consider its dual $B$ and let $V$ be the duality as in Equation \ref{eqn:5.1a} in the beginning of this section. We have the equality $\Delta(a)=V(a\ot 1)V^*$ valid in $AB\ot A$ where $AB$ is the Heisenberg algebra. Now assume that $\psi$ is a linear functional on the Heisenberg algebra satisfying the KMS property
\begin{equation*}
\psi(bx)=\psi(xS^2(b))
\end{equation*}
for all $b\in B$ and $x\in AB$. Then the restriction of $\psi$ to $A$ is right invariant.
\eprop

\bew
We use a Sweedler type notation for one of the copies of $V$ in the formulas below. For any $a$ in $A$ we find
\begin{align*}
(\psi\ot\iota)(\Delta(a))
&=(\psi\ot\iota)(V(a\ot 1)V^*)\\
&=(\psi\ot\iota)((v_{(1)}a\ot v_{(2)})V^*)\\
&=(\psi\ot\iota)(a\ot v_{(2)})V^*(S^2(v_{(1)})\ot 1)).
\end{align*}
Now we use that $V^*=(S\ot \iota)V$ (two times) and that $S$ is an anti-homomorphism. We get
\begin{align*}
(\psi\ot\iota)(\Delta(a))
&=(\psi\ot\iota)((a\ot 1)(S\ot\iota)((S(v_{(1)})\ot  v_{(2)})V))\\
&=(\psi\ot\iota)((a\ot 1)(S\ot\iota)(V^*V))=\psi(a)1.
\end{align*}
This completes the proof.
\ebew

Before we continue, let us make a couple of remarks about this result.

\opm
i) We have mentioned before that in finite dimensions, it is better to consider the non-involutive case. The above property and its proof work equally well for that  case. We just have to replace $V^*$ by $V\inv$. It is in fact the underlying idea of the existence proof of integrals on finite-dimensional Hopf algebras as given in \cite{VD4n}.\\
ii) One has to be a little more careful though. In principle it could happen that the restriction of $\psi$ to $A$ is trivially zero. This however can be overcome if we replace $\psi$ by $\psi(\,\cdot\,y)$ where $y$ is an element in the Heisenberg algebra that commutes with all elements of $B$. Such a functional will have the same KMS property. Remark that the Heisenberg algebra is in fact a full matrix algebra here. So it is no problem to find enough such elements.\\
iii) Assume that $\psi_0$ is a right integral on $A$. By the uniqueness of integrals, we will have a linear functional $\rho$ on the commutant of $B$ such that $\psi(ay)=\psi_0(a)\rho(y)$ for all $a\in A$ and $y\in B$. Therefore one can find a right integral for an element $y$ such that $\rho(y)\neq 0$.\\
iv) This property is intimately related with the result in Proposition 1.10 in \cite{Va-VD}.
\eopm

\prob\label{probl:5.22} 
Assume that $V$ is a multiplicative unitary. Assume that it has manageable antipodes. The $\sigma$-weak closures of the legs are von Neumann algebras $M$ and $\widehat M$. They carry coproducts. On $M$ we have $\Delta(x)=V(x\ot 1)V^*$. Formulate conditions so that it is possible to prove the existence of the Haar weights. Remark that it is sufficient to have either the left or the right Haar weights as we have the unitary antipode $R$ from the polar decomposition $S$ on $M$.
\eprob

%%%%%% Locally compact quantum groups
  
\section{\hspace{-17pt}. Locally compact quantum groups}\label{s:lcqgrps} % \input artikel6.tex

We begin the discussion in this section with formulating the definition of a locally compact quantum group in the setting of von Neumann algebras as it is given in the work of Kustermans and Vaes \cite{Ku-V3}.

\defin\label{defin:6.1}
Let $M$ be a von Neumann algebra and $\Delta$ a coproduct on $M$. The pair $(M,\Delta)$ is called a \emph{locally compact quantum group} if there  exists a left and a right Haar weight.
\edefin

A \emph{coproduct} $\Delta$ on the von Neumann algebra $M$ is by definition a normal unital $^*$-homomorphism from $M$ to the von Neumann algebra tensor product $M\ot M$, satisfying coassociativity $(\Delta\ot\iota)\Delta=(\iota\ot\Delta)\Delta$. A \emph{left Haar weight} is a faithful normal semi-finite weight on the von Neumann algebra that is left invariant and a \emph{right Haar weight} is a faithful normal semi-finite weight that is right invariant. A weight $\varphi$ is called \emph{left invariant} if $\varphi((\omega\ot\iota)\Delta(x))=\omega(1)\varphi(x)$ for all normal positive linear functionals $\omega$ on $M$ and all positive elements $x$ in $M$ with the property that $\varphi(x)<\infty$. Similarly for right invariance of a weight. We refer to \cite{Ku-V3} for details, see also \cite{VD9n}. 
\snl
The definition is found in \cite{Ku-V3} as well as in \cite{VD9n}. The advantage of the treatment in \cite{VD9n} over \cite{Ku-V3} is that in the first paper, the theory is developed \emph{without reference to the original C$^*$-approach}. We will also discuss this later in the present section.
\snl
The Haar weights are unique (up to a scalar) if they exist. For this reason, it is not necessary to include the symbols for the weights in the definition. This is what we have done also for Hopf algebras (cf.\ Definition \ref{defin:2.1} in Section \ref{s:fqgrps}) and multiplier Hopf algebras (cf.\ Definition \ref{defin:3.11} in Section \ref{s:dac}). 
\snl
For references on von Neumann algebras in general we refer to the item \emph{Basic references} in the introduction. For weights on von Neumann algebras and the related theory of Hilbert algebras, we refer to Chapters VI and VII in the book of Takesaki \cite{Tak2}.
\snl
We will explain further in this section why we start with the definition in the framework of von Neumann algebras and how to develop the theory from this. See Remark \ref{opm:6.3} and Remark \ref{opm:6.11a}. We first turn our attention to Kac algebras.

\subsection{A special case: Kac algebras}
Before we continue, let us compare this notion with  the one of a Kac algebra (see Definition 2.2.5 in \cite{E-S2}). We use a slightly different terminology together with notations compatible with the ones used elsewhere in this note.

\defin\label{defin:6.2}
A \emph{Kac algebra} is a pair $(M,\Delta)$ of a von Neuman algebra with a coproduct. It is assumed that there is \emph{a unitary antipode} $R$, defined here as an involutive $^*$-anti-automorphism that flips the coproduct. Furthermore there is a \emph{left Haar weight}.  It is defined as a faithful normal semi-finite weight $\varphi$ on the von Neumann algebra satisfying 
\begin{equation}
(\iota\ot\varphi)((1\ot y^*)\Delta(x))=R((\iota\ot\varphi)(\Delta(y^*)(1\ot x))\label{eqn:6.1a}
\end{equation}
for all $x,y\in M$ satisfying $\varphi(x^*x)<\infty$ and $\varphi(y^*y)<\infty$. For this formula to have a meaning, we need that $\varphi((\omega\ot\iota)\Delta(x))<\infty$ whenever $\omega$ is a positive element in $M_*$ and $x$ is a positive element in $M$ satisfying $\varphi(x)<\infty$. Finally it is assumed that $R\sigma^\varphi_t=\sigma^\varphi_{-t}R$ for all $t\in \mathbb R$ where $\sigma^\varphi$ is the modular automorphism group of $\varphi$.
\edefin

Let us first comment on the last condition.

\opm i) Because the unitary antipode $R$ flips the coproduct, the weight $\psi$, defined as the composition $\varphi\circ R$, is a right Haar weight when $\varphi$ is a left Haar weight. The modular automorphism groups are then related by the equation $R(\sigma^\varphi_t(x))=\sigma^\psi_{-t}(R(x))$ for all $t$ in $\mathbb R$ and all $x$ in the von Neumann algebra. Therefore, the last condition implies that the modular automorphism groups of the left and the right Haar weight coincide. \\
ii) The modular automorphism groups are also related by means of the modular element $\delta$ that appears as the Radon Nikodym derivative of the right Haar weight with respect to the left Haar weight. The relation is $\sigma^\psi_t(x)=\delta^{it}\sigma^\varphi_t(x)\delta^{-it}$ for all $t$ and all $x$.  So for a Kac algebra, the modular element is central.
\eopm

If $G$ is a locally compact group, both the function algebra $L^\infty(G)$ and the group von Neumann algebra VN($G)$ are Kac algebras. The last condition is obviously fulfilled in both cases. The function algebra is abelian so that the modular automorphism groups are trivial and the modular element is central. The group algebra is coabelian so that the left and the right Haar weights are the same. Hence their modular automorphism groups coincide and the modular element is trivial.

Also any finite quantum group is a Kac algebra because not only $S^2=\iota$, but also the left and right integrals coincide. In fact the integral is a trace so that the modular automorphism group is trivial. See Proposition \ref{prop:2.11}.
\snl
Any Kac algebra is a locally compact quantum group as in Definition \ref{defin:6.1} with the extra property that the antipode is a $^*$-map. The converse is probably not true because of the above remarks.

\prob\label{probl:6.8}
Find examples of locally compact quantum groups with an involutive antipode that are not Kac algebras.
\eprob

\subsection{Two equivalent approaches to locally compact quantum groups}
There is also a definition in the C$^*$-algebraic setting, but that is slightly more complicated. It is found in the original work of Kustermans and Vaes \cite{Ku-V1, Ku-V2}. 

\defin\label{defin:6.5b}
Let $A$ be a C$^*$-algebra and $\Delta$ a coproduct on $A$. Assume that  the linear span  of elements of the form
\begin{equation*}
(\omega\ot\iota)\Delta(a)
\tussenen
(\iota\ot\omega)\Delta(a),
\end{equation*}
where $a\in A$ and $\omega\in A^*$, are both dense subsets of $A$. If there exist a faithful left and a faithful right Haar weight, the pair $(A,\Delta)$ is called a \emph{locally compact quantum group}.
\edefin

Here a coproduct is a non-degenerate $^*$-homomorphism from $A$ to the multiplier algebra $M(A\ot A)$ of the minimal C$^*$-tensor product $A\ot A$. The Haar weights are lower semi-continuous densely defined KMS-weights. A left Haar weight is such a weight $\varphi$ that is left invariance. This means that 
\begin{equation*}
\varphi((\omega\ot\iota)\Delta(a))=\omega(1)\varphi(a)
\end{equation*}
for all positive elements $a$ in  $A$ such that $\varphi(a)<\infty$ and all positive $\omega$ in $A^*$.
\snl
Observe the presence of these density conditions as part of the axioms. They do not appear in the von Neumann algebraic definition. Similar properties are true also in that case, but they follow from the other axioms. The requirement of having  KMS-weights can be weakened, but in the end, it turns out that the Haar weights have this property.
\snl
There is a standard procedure to associate a C$^*$-algebraic locally compact quantum group to a von Neumann algebraic one and vice versa. The two procedures are inverses of each other in the sense that, if you apply them one after the other, you get back to the original locally compact quantum group. 
\snl
It is relatively easy to pass from a C$^*$-algebraic locally compact quantum group to a von Neumann algebraic one, without first developing the full theory in that setting. This has been argued in an appendix of \cite{VD9n}. One of the main results achieving this is Proposition A.6 of that appendix. The indications given for the proof of that proposition are not completely correct. It is the intention to give a better and more complete argument in a newer version of these notes \cite{VD11n}. This is related with Problem \ref{probl:6.10a} formulated further.
\snl
The other direction is more involved. Given a von Neumann algebraic locally compact quantum group $(M,\Delta)$ as in Definition \ref{defin:6.1}, it is relatively easy to get  the associated C$^*$-algebra $A$. It is the norm closure of the right leg of the multiplicative unitary $V$ constructed from the right Haar weight on $M$. It requires some work to show that this norm closure is self-adjoint. The antipodes must be constructed and one has to show that they are manageable as in Definition \ref{defin:5.11}. Then one can use Proposition \ref{prop:5.12a}. The next step is to show that the restriction of the coproduct to $A$ is indeed a coproduct on the C$^*$-algebra. Finally it is needed that the restrictions to $A$ of the original Haar weights on the von Neumann algebra $M$ are still densely defined. Some aspects are standard but not all of them. The results are found in \cite{Ku-V3}, see also a remark in Appendix A of \cite{VD9n}. In \cite{VD11n} we plan to cover this aspect in detail. 
 \snl
Because of this, it is correct to say that the two approaches to locally compact quantum groups yield the same objects. 

\opm\label{opm:6.3}
i) The development of the theory is \emph{easier} in the von Neumann algebraic framework. There are several reasons for this. We collect them here.
\begin{itemize}
\item[-] The notion of a coproduct on a von Neumann algebra is easier than the one on a C$ ^*$-algebra.
\item[-] The theory of weights on von Neumann algebras is easier and better documented in the literature.
\item[-] The set of axioms in Definition \ref{defin:6.1} is simpler than in Definition \ref{defin:6.5b}.
\item[-] When using von Neumann algebras, one can use Hilbert space techniques.
\end{itemize}
ii) On the other hand, it is \emph{more natural} to formulate the notion with C$^*$-algebras. Indeed, we think of a C$^*$-algebra as a locally compact quantum space while we use von Neumann algebras for quantized measure spaces. 
\eopm
The earlier theory of Kac algebras is a von Neumann algebraic one. This inspired Masuda and Nakagami to develop a theory of locally compact quantum groups with von Neumann algebras in \cite{M-N} (1994). On the other hand, the quantum $SU_q(2)$ as introduced by Woronowicz in \cite{W1} (1987), is formulated in C$^*$-algebraic terms. The same is true for the compact quantum groups in \cite{W2,W3}.  Following the consensus here, the search for a good theory of locally compact quantum groups was further done in the C$^*$-algebraic setting. This led to the work of Masuda, Nakagami and Woronowicz \cite{M-N-W} that appeared in 2003 (but with preliminary versions available many years before that) and simultaneously the work of Kustermans and Vaes \cite{Ku-V1, Ku-V2} (in 1999 and 2000). 
\snl
Later, Kustermans and Vaes also treated their locally compact quantum groups in the von Neumann algebraic framework (\cite{Ku-V3}, 2003). The approach heavily depends on their original paper \cite{Ku-V2} were the theory is developed in the C$^*$-algebraic setting. This is different in \cite{VD9n} where the von Neumann algebraic approach is obtained without reference to the C$^*$-algebraic one of \cite{Ku-V2}. 
\snl
Before we continue, we need to say a few more words about the comparison of the work by Masuda, Nakagami and Woronowicz \cite{M-N, M-N-W} on the one hand and that of Kustermans and Vaes on the other hand \cite{Ku-V2,Ku-V3} on the other hand. 

\opm
The main difference between the two approaches is that in the first work, the antipode with its polar decomposition is assumed while in the second one, it is proven to exists from the axioms. To do this, it is necessary to assume both the existence of a left and a right Haar weight while in the first work, only a left Haar weight is used. There the right Haar weight is obtained by composing the left Haar weight with the unitary antipode that appears in the polar decomposition of the antipode. 
\eopm

It is now generally accepted that the theory of Kustermans and Vaes is stronger. Indeed, the right Haar weight is easy to obtain from the existence of the antipode with its polar decomposition. % {\blauw Dit statement goed nakijken}.

\subsection{The theory of  locally compact quantum groups}
We will finally discuss some aspects of the treatment of locally compact quantum groups in the von Neumann algebraic framework. In particular, we will show how the knowledge of the algebraic quantum groups, recalled in Section \ref{s:alg}, inspires this treatment. We also will formulate some problems on the way. 
\snl
We will not cover other aspects of the theory here. They are found in the literature. Moreover, once the antipode is obtained, much of the other steps are now more or less standard. The approach in \cite{VD9n} is most adapted to the discussion here. We also plan a new expanded version of these notes (cf.\ \cite{VD11n}). 
\snl
The starting point is a von Neumann algebraic quantum group as in Definition \ref{defin:6.1}. We denote a left Haar weight with $\varphi$ and a right Haar weight with $\psi$. 
\snl
As should be clear from the remarks made earlier, the crucial step is the construction of the antipode. In our approach, the underlying idea is taken from Proposition \ref{prop:4.10a}, characterizing the antipode in a multiplier Hopf algebra. We reformulate it for a Hopf $^*$-algebra $(A,\Delta)$ where the property reads as follows. 

\prop
If $a\in A$ we can write $a\ot 1 =\sum_i \Delta(p_i)(1\ot q_i^*)$ and then we have $$S(a)^ *\ot 1=\sum_i \Delta(q_i)(1\ot p_i^*).$$
\eprop

It is possible to formulate the result for multiplier Hopf $^*$-algebras if we allow some kind of approximation of $a\ot 1$ by elements of the type $\sum_i \Delta(p_i)(1\ot q_i^*)$. Then we get an approximation of $S(a)^*\ot 1$. 
\snl
We can illustrate the same phenomenon in the case of a locally compact group $G$ where we take  the $^*$-algebra $C_c(G)$ of continuous complex functions with compact support. Indeed, if $f\in C_c(G)$ we can approximate
\begin{equation*}
f(r)=f(rs\cdot s\inv)\simeq\sum_ip_i(rs)\overline{q_i(s\inv)},
\end{equation*}
and we get 
\begin{equation*}
\overline{f(r\inv)}=\overline{f(s\cdot(rs)\inv)}=\sum_iq_i(rs)\overline{p_i(s)}.
\end{equation*}

This idea is now used in \cite{VD9n} to obtain the following candidate for the antipode on the von Neumann algebra $M$ (see Definition 1.2 in \cite{VD9n}). We assume that the von Neumann algebra acts on the Hilbert space $\mathcal H$.

\defin\label{defin:6.6a}
For an element $x\in M$ we say that $x\in\mathcal D$ if there is an element $x_1 \in M$ satisfying the following condition. For all $\varepsilon > 0$ and vectors $\xi_1, \xi_2, \dots, \xi_n,\eta_1, \eta_2, \dots, \eta_n$ in $\mathcal H$, there exist elements $p_1, p_2, \dots, p_m, q_1, q_2, \dots, q_m$ in $M$ such that
\begin{align} 
&\|x \xi_k\ot \eta_k - \sum_j \Delta(p_j)(\xi_k\ot q_j^*\eta_k)\| < \varepsilon \\
&\|x_1 \xi_k\ot \eta_k - \sum_j \Delta(q_j)(\xi_k\ot p_j^*\eta_k)\| < \varepsilon
\end{align}
for all $k$.
\edefin

In order to use this to obtain the antipode we first need  to prove that the element $x_1$ is unique if it exists. Then we can define $S$ on $\mathcal D$ by $S(x)^*=x_1$. Next we must show there exists elements $x$ with this property and that in fact the set $\mathcal D$ is dense in $M$. 
\snl
We know from the algebraic theory discussed in Section \ref{s:alg} that we expect to need the right Haar weight for the first problem. Recall that the canonical map $p\ot q\mapsto \Delta(p)(1\ot q)$ is shown to be injective using the faithful right integral. We need a similar result here to get uniqueness of the element $x_1$ if it exists. On the other hand, we also know from the algebraic theory that it is expected to get enough elements in $\mathcal D$ by the existence of the left integral. See the discussion in Subsection \ref{ss:LST} of Section \ref{s:alg} on the Larson Sweedler theorem.
\snl
It is not hard to imagine that the case here is far more complicated than in the algebraic setting, mainly due to the fact that we are working with approximations. 
\snl
It can be done to define the antipode in this way, but it turns out to be better to first define the map $x\mapsto S(x)^*$ on the Hilbert space level. 
\snl
The idea is as follows. Consider the GNS Hilbert space $\mathcal H_\psi$ constructed from the right Haar weight. The canonical map, defined from the subset $\mathcal N_\psi$ of elements $x$ in $M$ satisfying $\psi(x^*x)<\infty$ to the Hilbert space is denoted by $\Lambda_\psi$.  We now define the canonical map $T_1$ on the Hilbert space level with the formula inspired by the one in Proposition \ref{prop:4.9} of Section \ref{s:alg} in the case of an algebraic quantum group. It yields an isometric linear map $V$ from the Hilbert space tensor product $\mathcal H_\psi\ot \mathcal H_\psi$ to itself. The map $V$ turns out to be unitary, but this is not immediately clear. We refer to Problem \ref{probl:6.10a}. Still one can now define an operator $K$ with domain $\mathcal D(K)$ on $\mathcal H_\psi$, along the lines of Definition \ref{defin:6.6a}.

\defin
Let $\xi\in \mathcal H_\psi$. We say that $\xi\in \mathcal D(K)$ if there is a vector $\xi_1\in \mathcal H_\psi$ satisfying the following condition: For all $\varepsilon > 0$ and vectors $\eta_1, \eta_2, \dots, \eta_n$ in $\mathcal H_\psi$, there exist elements $p_1, p_2, \dots, p_m, q_1, q_2, \dots, q_m$ in $\mathcal N_\psi$ such that
\begin{align} 
&\|\xi\ot \eta_k - V(\sum_j \Lambda_\psi(p_j) \ot q_j^*\eta_k)\| < \varepsilon \\
&\| \xi_1\ot \eta_k - V(\sum_j \Lambda_\psi(q_j)\ot p_j^*\eta_k)\| < \varepsilon
\end{align}
for all $k$.
\edefin

It is now shown that the vector $\xi_1$ is unique if it exists in Lemma 1.7 of \cite{VD9n} by techniques well-established in the theory of normal weights on von Neumann algebra. Just observe that, as expected, the right Haar weight is used here. 
\snl
This allows to define the linear operator $K$ on $\mathcal D(K)$ by $K\xi=\xi_1$. It is immediately clear from the definition that it is a closed conjugate linear involutive operator. For the density of the domain $\mathcal D(K)$, again as expected, the left Haar weight is used. See Lemma 1.12 of \cite{VD8n} and compare the proof of that lemma with arguments used in Subsection \ref{ss:LST} on the Larson Sweedler theorem. Remark in passing that the density of the domain of $K$ is proven simultaneously with the argument that the operator $V$ is not just isometric, but unitary.
\snl
Before we continue, we formulate the following related problem.

\prob\label{probl:6.10a}
Given the right invariant Haar weight $\psi$ on the pair $(M,\Delta)$, it uses standard techniques to define the canonical map $V$ on the Hilbert space $\mathcal H_\psi\ot\mathcal H_\psi$ and to show that it is an isometry. However, it is not trivial to show that it is actually a unitary. This is easy in the algebraic setting using the surjectivity of the canonical map $a\ot a'\mapsto \Delta(a)(1\ot a')$. But it is not so in the topological setting. In the original  work of Kustermans and Vaes \cite{Ku-V2} it is `Kusterman's trick' that does the job (see the remarks above). This is also how it is done in \cite{VD9n}. It would be nice to find a more direct and more natural way to prove this property. A good solution here is also important for the passage from C$^*$-algebras to von Neumann algebra (see \cite{VD11n}). At a workshop in Marseille some years ago, an elegant solution to this problem has been presented by P.\ Kasprzak, but there are no traces of his proof unfortunately.
\eprob

Once it is shown that the operator $K$ is well and densely defined, we can proceed to complete the construction of the antipode. It is in fact done by \emph{redefining} the domain of $S$. We say that $x\in\mathcal D(S)$ if $\xi\in\mathcal D(K)$ implies that $x\xi\in\mathcal D(K)$ and $K(x\xi)=S(x)^*K\xi$. The idea behind this definition is that $S$ is an anti-homomorphism so that $S(xy)^*=S(x)^*S(y)^*$ and if we look at this equation on the Hilbert space level where $K$ stands for $x\mapsto S(x)^*$, this formula reads as $K(x\xi)=S(x)^*K(\xi)$.  It is not clear if this new definition of the domain $\mathcal D(S)$ is the same as the one $\mathcal D$, defined earlier in Definition \ref{defin:6.6a} but this is not important for the further development of the antipode and the remaining of the theory.
\snl
By defining the antipode via the operator $K$, we can take advantage of the polar decomposition of the closed linear operator $K$ to get the polar decomposition of the antipode. Further Hilbert space methods are used to proceed.
\snl
For details about the construction of the antipode and how to obtain the polar decomposition of it, we refer to Section 1 of \cite{VD9n}.

This is the core of the further development of the theory.

\opm\label{opm:6.11a}
In Remark \ref{opm:6.3} we formulated several reasons for the fact that the von Neumann algebraic theory is easier than the C$^*$-algebraic one. One of them was the possibility to use Hilbert space techniques. This is indeed a great advantage for treating the antipode further. It is done via the polar decomposition of the closed operator $K$.
\eopm
It is however not completely correct to use this argument because Hilbert space techniques and the polar decomposition of the operator $K$ are also used in the C$^*$-algebraic development. But still, they are more natural in the development of the von Neumann algebraic theory.
\snl
It would take us now too far to comment more on the further development of the theory of locally compact quantum groups. But we have made our point. We have shown how the algebraic theory is a source of inspiration to develop the analytical one. This is illustrated with the construction of the antipode above. And as it is well-known, this step is the crucial one making further results with their proofs possible. We refer again to \cite{VD9n} en \cite{VD11n}.

%%%%%% Conclusions
  
\section{\hspace{-17pt}. Conclusions} \label{s:conclusions} % \input artikel7.tex

In this paper we discussed the passage from the theory of Hopf algebras to  locally compact quantum groups. This starts with finite quantum groups. They are defined as finite-dimensional Hopf $^*$-algebras with an operator algebra as the underlying algebra. It then moves via discrete and compact quantum groups to algebraic quantum groups. Algebraic quantum groups are defined as multiplier Hopf $^*$-algebras with positive integrals. Discrete and compact quantum groups are in the first place considered as special cases of algebraic quantum groups. The positivity of the integrals allows a Hilbert space realization using multiplicative unitaries. These are first studied from a purely algebraic point of view to clarify the link between the algebraic theory and the operator algebraic one. Finally we arrive at locally compact quantum groups, both in the C$^*$-algebraic setting as well as in the von Neumann algebra framework.
\snl
Along the way, we have indicated what kind of choices have to be made and how this has led to different and sometimes annoying conventions with different notations and formulas. We have included some historical information about the development of the theory. 
\snl
Also some remaining problems have been mentioned. Among them, the most challenging one  is the existence problem of the Haar weights. There are satisfactory theorems for finite, discrete and compact quantum groups, but not beyond that. It seems to be a difficult problem and we have provided some ideas in Section \ref{s:mu} on multiplicative unitaries (see Problem \ref{probl:5.22}). Another problem of interest is to find examples of locally compact quantum groups with involutive antipode that are not Kac algebras (see Problem \ref{probl:6.8}).
\snl
It is my opinion that a comprehensive textbook or lecture notes on the following two topics would be most welcome. 
\snl
First there is a need for a complete work on the theory of algebraic quantum groups, as understood in this paper. It should include all features, earlier obtained for (finite-dimensional) Hopf algebras and also true for multiplier Hopf algebras and algebraic quantum groups. At this moment, many results are spread over various papers, dealing with multiplier Hopf algebra, algebraic quantum groups and beyond. Moreover these more general theories sometimes gave new insights of known aspects of the original papers and should be included as such in the development.
\snl
There is a need for a work on locally compact quantum groups as well. 
The writing of it should be used to settle the terminology and conventions and where needed, with comparison  to the literature when different conventions are used. 
\snl
There exist articles of this type. We have a survey paper on multiplier Hopf algebras (see  \cite{VD-Z}, 2000), but that certainly should be updated with more recent results and new approaches on earlier results. There is also the book of Timmermann (see \cite{Timmermann}, 2008). This also contains a survey on multiplier Hopf algebras and algebraic quantum groups and duality but this is very concise and contains almost no proofs. Moreover since then new research has been done. 
Finally, I believe it is desirable for analysts working with the operator algebra approach to quantum groups to  keep in contact with the algebraists in the field of Hopf algebras. And vice versa. We can learn from each other.
\nl\nl

\end{document}